\documentclass[12pt,a4paper]{article}

\usepackage{pdflscape}
\usepackage{amsmath}
\usepackage{amssymb}
\usepackage{latexsym}
\usepackage{srcltx}
\usepackage{graphics}
\usepackage{color}
\usepackage{epsfig}
\usepackage{color}
\usepackage[ruled,vlined]{algorithm2e}

\usepackage{subcaption}

\usepackage{anyfontsize}
\usepackage{tabularx}

\usepackage{graphicx}
\usepackage{multirow}
\usepackage{array}

\usepackage{amssymb, amsmath, amsfonts, mathtools}
\usepackage{esint}
\usepackage{graphicx}
\usepackage{breqn}
\usepackage{float}
\usepackage{caption}
\usepackage{subcaption}
\usepackage{hhline}
\usepackage{multirow}
\usepackage{rotating}
\usepackage{fixltx2e}

\newcommand{\co}{{\mathbb C}}

\newcommand{\re}{{\mathbb R}}
\newcommand{\rat}{{\mathbb Q}}
\newcommand{\n}{{\mathbb N}}

\newcommand{\z}{{\mathbb Z}}
\newcommand{\cA}{{\mathcal{A}}}

\newcommand{\cV}{{\mathcal{V}}}

\newcommand{\cT}{{\mathcal{T}}}
\newcommand{\cH}{{\mathcal{H}}}

\newcommand{\cR}{{\mathcal{R}}}

\newcommand{\cP}{{\mathcal{P}}}

\newcommand{\by}{{\boldsymbol y}}
\newcommand{\bx}{{\boldsymbol x}}
\newcommand{\nul}{{\boldsymbol 0}}

\newcommand{\ba}{{\boldsymbol a}}
\newcommand{\balpha}{{\boldsymbol \alpha}}
\newcommand{\bb}{{\boldsymbol b}}
\newcommand{\bc}{{\boldsymbol c}}
\newcommand{\bz}{{\boldsymbol z}}
\newcommand{\be}{{\boldsymbol e}}

\newcommand{\bv}{{\boldsymbol v}}

\newcommand{\bC}{{\boldsymbol{C}}}

\newtheorem{theorem}{Theorem}
\newtheorem{prop}{Proposition}
\newtheorem{lemma}{Lemma}
\newtheorem{cor}{Corollary}
\newtheorem{remark}{Remark}
\newtheorem{ex}{Example}
\newtheorem{defi}{Definition}

\date{}

\author{
Vladimir Yu.~Protasov
\thanks{University of L'Aquila, Italy,  {e-mail: \tt\small
vladimir.protasov@univaq.it}}}

\title{The Barabanov norm is generically unique, \\ 
simple, and easily computed\thanks{
This work is supported by the RFBR  grants no 19-04-01227 and 20-01-00469
}
}

\begin{document}
\maketitle

\begin{abstract}

Every irreducible discrete-time linear switching system possesses an 
invariant convex Lyapunov function (Barabanov norm), which  provides 
a very refined analysis of trajectories. 
Until recently that notion remained rather theoretical apart from special cases. 
In 2015 N.Guglielmi and M.Zennaro showed that many systems possess at least one 
simple Barabanov norm, which moreover, can be efficiently computed. In this paper we 
classify all possible Barabanov norms for discrete-time systems. We 
prove that, under mild assumptions, such norms are unique and are either 
piecewise-linear or piecewise quadratic. Those assumptions can be verified 
algorithmically and the numerical experiments show that a vast majority of systems 
satisfy them. For some narrow classes of systems,  
there are more complicated Barabanov norms  but they can still be classified and constructed. 
Using those results we find all trajectories of the fastest growth. 
They turn out to be eventually periodic with special periods. Examples and numerical results are presented.

\bigskip

\smallskip

\begin{flushright}
\noindent  \textbf{AMS 2010} {\em subject
classification:  93D20, 39A22, 46B20, 52B12}
\end{flushright}

\end{abstract}
\bigskip

\begin{center}
\large{\textbf{1. Introduction}}
\end{center}
\bigskip

Discrete-time linear switching systems of the form 
\begin{equation}\label{eq.sys}
\left\{
\begin{array}{l}
\bx(k+1)  \ = \ A(k)\bx(k);\\
  A(k) \in \cA \, , \quad  k \in \z_+\, ; \\
\bx(0) \ =  \ \bx_0 
\end{array}
\right. 
\end{equation} 
are the subject of an extensive literature. The set $\cA$ is a compact set of $d\times d$ matrices. In this paper we deal with finite families $\cA = \{A_1,\ldots , A_m\}, \, m \ge 1$, 
and often identify the system with the corresponding family. 
The sequence of matrices $\{A(k)\}_{k=0}^{\infty}$ taken from~$\cA$ 
(with repetitions permitted) is the {\em switching law}. 
The sequence of points $\{\bx(k)\}_{k=0}^{\infty}$
from~$\re^d$ satisfying~(\ref{eq.sys}) for some switching law~$A(\cdot)$ is a {\em trajectory} of the system. The switching law,  
along with the initial point~$\bx_0$,  defines the trajectory. 

 The fastest possible growth of trajectories as $k \to \infty$ is an important issue in many applied problems. It is closely related to the asymptotic stability of linear and non-linear systems~\cite{Gurv, L,  MP2},  to the regularity exponents of fractal curves and surfaces~\cite{CHM, P06},  wavelets, and subdivision schemes~\cite{G, PG15}, 
 to the growth of special sequences in combinatorics, number theory, and the theory of formal 
 languages, in the automata theory, etc.,  see~\cite{BCJ, JPB, MOS, P17} and references therein. If the system is {\em irreducible}, i.e., 
 the matrices~$A_1, \ldots , A_m$ do not share a nontrivial invariant 
linear  subspace, then the maximal value of $\|\bx(k)\|$ over all trajectories 
(with fixed~$\bx_0$) 
 is asymptotically equivalent to~$\rho^{\,k}$. More precisely, it is 
 between~$C_1 \|\bx_0\|\, \rho^{\,k}\, $ and $\, C_2\|\bx_0\| \, \rho^{\,k}$, where 
 $C_1 \le C_2$ are positive constants and $\rho = \rho(\cA)$ is the {\em joint spectral radius} (JSR) of the family~$\cA$.  We recall the definition of the JSR below. There are efficient methods to estimate the JSR~\cite{AS, G, J} and, in many cases, even to compute it precisely~\cite{GP1}. So, the exponent of the fastest growth~$\rho$ can be computed. However, 
 it is not enough to have a comprehensive information on the growth of trajectories, 
 since $C_1$ can be very small or $C_2$ very large. One needs to estimate the constants $C_1$
 and $C_2$.  
 This problem, however, is more difficult. Even their rough estimations are usually hard.  Theoretically this problem can be solved by 
 using the {\em invariant convex Lyapunov function} also called the  
 {\em Barabanov norm}. 
 \begin{defi}\label{d.bar}
 An invariant convex Lyapunov function (Barabanov norm) of a family 
 of matrices $\cA = \{A_1, \ldots , A_m\}$ is a norm $f$ in $\re^d$ such that 
 \begin{equation}\label{eq.bar}
 \max_{i=1, \ldots , m}f(A_i\bx)\ = \ \rho\, f(\bx)\quad  \mbox{for every} \ \bx \in \re^d, 
 \end{equation} 
 where $\rho = \rho(\cA)$ is the joint spectral radius of~$\cA$. 
 \end{defi}
 By iterating equation~(\ref{eq.bar}) we obtain $\max_{s_1, \ldots , s_k}f(A_{s_k}\cdots A_{s_1}\bx)\, = \, \rho^{\,k}\, f(\bx)$ for every~$k$. Consequently, for the Barabanov norm, we have $C_1 = C_2 = 1$, which means that this norm is optimal among all possible norms in~$\re^d$. For an arbitrary norm~$\|\cdot \|$, say, Euclidean,
the constants~$C_1, C_2$ can be obtained by the maximal and minimal values of $f(\bx)/\|\bx\|$, provided the Barabanov norm $f$ is known. Thus, if the Barabanov norm is available, 
then the problem of estimating the maximal growth of trajectories in every norm is 
efficiently solved.   
Moreover, in that case  it is possible to find all switching 
laws $A(\cdot)$ corresponding to the fastest growth of trajectories, 
 for which $\limsup_{k \to \infty}\rho^{-k}\|A(k)\ldots A(1)\| \, > \, 0$.  

Thus, in the analysis of the trajectories, it is very desirable to have a Barabanov norm.  
It was shown in~\cite{B} that such a norm does exist for every irreducible family~$\cA$. 
This is a purely existence result, all of its known proofs are non-constructive. 
Is it possible 
to obtain the Barabanov norm in a closed form? There are several arguments saying that the answer 
should be negative:  
\smallskip 

1) {\em The non-uniqueness.} We always consider  the uniqueness of norms 
up to their multiplication by a constant. Simple examples in 
$\re^2$ already show that the Barabanov norm may not be unique. Say, if $\cA$ consists of one  $2\times 2$ matrix that defines a rotation of the plane 
by the right angle, then the $L_2$-norm  and the $L_1$-norm in~$\re^2$ are both Barabanov.    
 \smallskip 

2) {\em Non-convergence of the power sequence $F^{n}[g]$ as $n\to \infty$}, where 
$F$ is a map on the set of norms:~$\, F[g](\bx)\, = \, 
\rho^{-1}\max_{i=1, \ldots , m}g(A_i\bx)$. The Barabanov norm~$g = f$ is a fixed point for~$F$.  
However, the iterations $F^{n}[g]$ may not converge for some~$g$ even if that fixed point is unique.  For example, if $\cA$ consists of one rotation 
of~$\re^2$ by an angle $\alpha$ such that $\alpha/\pi$ is irrational, 
then there is a unique Barabanov norm~$f$, which is the Euclidean norm. 
However, for the $L_1$-norm  $g(x_1, x_2) = |x_1| + |x_2|$, the sequence 
$F^{n}[g]$ does not have a limit. This shows that even in case of uniqueness, the 
Barabanov norm cannot be computed by the power method. 
  \smallskip 

3) {\em The fractal-like boundary of the unit sphere.} Let 
$B = \{\bx \in \re^d \ | \ f(\bx) \le 1\}$ be a unit ball of the 
Barabanov norm and let $G = B' = \{\by \in \re^d \ | \ \max_{\bx \in B} (\bx, \by) \le 1\}$
 be its polar. Then, as it was proved in~\cite[Theorem 17]{PW},  the convex hull 
 of images~$A_i^TG, i = 1, \ldots , m$, is homothetic to~$G$ itself: 
\begin{equation}\label{eq.invarT}
 {\rm co}\, \Bigl(\,  \bigcup_{i=1}^{m} \, A_i^TG\, \Bigr)\ = \ \rho\, G, 
\end{equation} 
 where $A_i^T$ denotes the transpose  matrix to~$A_i$. 
 The convex body~$G$ possessing this property generates the {\em Protasov norm} 
 according to the terminology in~\cite{PW}. Thus, Barabanov's and Protasov's norms are 
 dual to each other. The  property~(\ref{eq.invarT})
 is closely related to the definition of self-similar fractals by J.Hutchinson~\cite{H}. 
 In fact, if we iteratively construct such a convex body~$G$ on the plane, 
 we will see that there should be segments on its boundary (because of taking the 
 convex hull) and those segments multiply with iterations. Hence, the boundary 
 should have a structure somewhat similar to the Cantor set. See~\cite{S} for more details. 
    \smallskip 
 
 The  arguments above suggest  that the notion of the Barabanov norm is rather theoretical 
 and can hardly be evaluated for general matrices. This is indeed a common belief among specialists 
 working in discrete-time switching systems. Many interesting theoretical results on 
 Barabanov's 
norms can be found in~\cite{K2, K4, LX, Ma08, M1, M2, PW}.   However, as it was remarked in 2012 by 
 R.\,Teichner and M.\,Margaliot: ``Although the Barabanov norm was studied extensively, it seems that there are only few examples where it was actually
computed in closed form''~\cite{TM}. 

Surprisingly, in 2015 N.\,Guglielmi and M.\,Zennaro~\cite{GZ} showed 
that for many maxtrix families (in particular, for all known families from applications) 
it is possible to construct at least one  Barabanov norm in an explicit form. 
This form is either piecewise linear or piecewise-quadratic and   it  can be 
found within finite time. For constructing that norm they put to good use  the {\em invariant polytope algorithm} from~\cite{GP1}, whose idea also traces back to works of 1996~\cite{P96} and of~2005~\cite{GWZ}. 
An earlier versions of that algorithm appeared in~\cite{BJP, GZ08}, see also~\cite{Ma08, P96} for 
other related algorithms. 

The invariant polytope  algorithm produces an {\em invariant 
convex body}~$G$ possessing the property
\begin{equation}\label{eq.invar}
 {\rm co}\, \Bigl( \, \bigcup_{i=1}^{m} \, A_iG\, \Bigr)\ = \ \rho\, G, 
\end{equation} 
and thus finds precisely the value of JSR $\rho(\cA)$. Actually the 
algorithm finds a {\em dominant product}~$\Pi = A_{s_n}\ldots A_{s_1}$
(see Definition~\ref{d.10} in Section~2) such that $|\lambda|^{1/n} = \rho(\cA)$, where 
$\lambda$ is the {\em leading}, i.e., the largest in modulus 
 eigenvalue of~$\Pi$. 
It was proved in~\cite{GP1} 
that the algorithm halts within finite time if and only if the product~$\Pi$ is 
dominant and its leading eigenvalue~$\lambda$ is unique and simple. 
In this case the obtained invariant body~$G$ is either a polytope (if $\lambda \in \re$)
or a convex hull of several ellipses (if $\lambda \notin \re$). 
Then from a result of E.Plischke and F.Wirth~\cite[Theorem 17]{PW}  it follows that if we get an invariant body~$G^*$ of the 
transpose family $\cA^* \, = \, \{A_1^T, \ldots , A_m^T\}$, then 
the function~$f(\bx) \, = \, \max_{\by \in G*}(\by, \bx)$ is a Barabanov norm. 

A lot of numerical experiments done in~\cite{GP1, GP2, M} demonstrate that, 
for a vast majority  of matrix families, the invariant polytope algorithm 
halts and hence produces an invariant body.  
There are well-known counterexamples~\cite{BTV, Sid1} but they are absolutely rare in practice. 
Having applied that algorithm to the transpose family we obtain 
the Barabanov norm. 

\smallskip 

\textbf{An assumption based on numerical experiments}. 
Let us clarify our claim on the  ``vast majority of matrix families
in the numerical experiments''. First of all, the invariant polytope 
algorithm is robust: if it halts for some family of matrices, than it does 
for all close families performing the same number of iterations. 
The parameters of robustness are efficiently estimated~\cite[Section~2.5]{GP1}. 
This makes it possible to avoid using exact arithmetics  or rational matrices 
in the numerical computations. All the experiments are performed
with  rounding using well-defined tolerance parameters. There are several computer implementations applying various software~\cite{M, PZ1}.  Several hundreds of numerical tests have been done  in 
 dimensions up to~$20$ with two sorts of matrix families: 1) families from known applications 2) randomly generated matrices. In all these experiments (100~\%) the invariant polytope 
 algorithm terminates within finite (usually quite short) time.
 See~\cite{GP1, GP2, M, PZ1} for more details.   The statistics of a small part of those 
 experiments is demonstrated in Section~10.  This allows us to assume that 
 a generic family of matrices possesses this property. 
 By ``generic'' we mean that for every $m$ and $d$, 
 the property holds for  an open set of full Lebesgue measure in the
 space~$\re^{md^2}$ (the space  of families of~$m$ matrices $d\times d$).  
 We are not aware of any rigorous  results approving this claim
 and we believe this is a challenging theoretical problem. Therefore, we make the assumption that for a generic family the algorithm halts based on numerical experiments. 
 For an arbitrary matrix family, this can be checked directly by running the algorithm.

\smallskip 

\textbf{Statements of the problems and a summary of main results}. 
Thus,  a generic matrix family has at least one Barabanov 
norm that can be found in a closed form. Its unit ball is either a polytope (if the leading eigenvalue~$\lambda$ of the 
dominant product is real) or a polar to a convex  hull of ellipses (otherwise). 
A question arises if it has other Barabanov norms and, if so, how many   and of what structure? 
If there are norms with fractal properties,  how to find them? 
And how to find an optimal one among all Barabanov norms? 
In this paper we answer all those questions. We prove that in most cases 
the Barabanov norm is unique. This means that for all 
generic families, Barabanov norms are simple (either piecewise-linear or 
piecewise-quadratic) and there are no others. They are easily computed as 
maxima of several linear (respectively, quadratic) functionals. Thus, the invariant polytope algorithm produces 
not some norm but all possible Barabanov norms. In particular, there are no 
``fractal-like'' norms among them. More precisely, we prove that the uniqueness takes place 
if the leading eigenvalue $\lambda$ 
is either real or complex with an argument~$\pi q$  with  irrational~$q$
(Theorems~\ref{th.10},~\ref{th.12} in Section~2). 
But what about the known simple examples when the Barabanov norm is not unique? 
It turns out that all of them belong to the third case: $\lambda$ 
is complex with an argument~$\pi q$, where~$q$ is a rational non-integer number. 
In this case, as we shall see, a family~$\cA$ has  a rich variety if Barabanov norms. 
Nevertheless, they all can be classified. We do it in Theorem~\ref{th.40}, Section~8.  Next we extend those results to 
families with several (more than one) dominant products. This case is important in applications
(see Section~6 for details). We prove that in this case there always exist 
infinitely many Barabanov norms but they are all quite simple and can be found by
a modified version of the invariant polytope algorithm. In Section~7 
we apply our results to the classification  of  trajectories  of the 
fastest growth. All of them  can be explicitly found: 
a switching law provides the fastest growth, i.e.,  $\|A(k)\cdots A(1)\| \, \ge \,  
C\rho^k, \, k \in \n$,  precisely when it is {\em eventually} periodic, 
 i.e., $A(k+n) = A(k)$ for all~$k>N$, where $n, N$ are some natural numbers, 
  and the period is equal to 
one of the dominant products. For all other trajectories, we have 
$\|A(k)\cdots A(1)\|\, \rho^{-k} \to 0$ as $k \to \infty$
 (Theorem~\ref{th.120}, Section~8). In Section~9 
we turn to positive systems, when the Barabanov norm is always unique and 
is piecewise-linear (provided a dominant product exists). A modification of the invariant polytope algorithm 
for positive systems is very efficient: it constructs the Barabanov norm 
even for very large dimensions~$d$ (several thousands). 
Finally, in Section~10 we present numerical results and discuss the computational issue. 
We will see that in most cases the time of constructing the Barabanov norm does not exceed 
that for constructing other Lyapunov functions by algorithms known from the literature.  
\smallskip 

\textbf{Novelty}. Our results can be divided into four main groups:
\smallskip 

1) The proof of uniqueness of the Barabanov norm provided  
the dominant product has a leading eigenvalue which is either real or complex 
with an irrational  $\, {\rm mod}\, \pi\, $  argument. This shows  that for a generic family of matrices, the Barabanov norm is unique, has a simple structure, and can be efficiently found (Sections~2 - 4).
In the remaining case (the non-real eigenvalue with a rational ${\rm mod}\, \pi$ argument), the family of matrices has a large variety of Barabanov norms. 
We classify them and present an algorithm to find them all (Section~8). 
  \smallskip 

2) In case of finitely many dominant products, we show that there always exists 
an infinite set of Barabanov norms but  all of them have a simple structure and can be explicitly found. An algorithm for their construction  
is presented (Section~6).

\smallskip 

3) For every discrete-time system with finitely many dominant products, 
all switching laws corresponding to the fastest growth of trajectories are explicitly found. 
All  trajectories of the fastest growth are classified 
(Section~7). 
\smallskip 

4) For positive systems (Section~9), we introduce the monotone Barabanov
and prove that they are  unique and piecewise-linear. They are
found by a modification of the invariant polytope algorithm, which works 
efficiently even in very large dimensions.
\smallskip

\textbf{Auxiliary facts and notation}. We use bold letters for vectors and standard 
letters for numbers, so $\bx = (x_1, \ldots , x_d)^T \in \re^d$. We consider a discrete-time system~(\ref{eq.sys})
in $\re^d$ with a finite family of matrices~$\cA = \{A_1, \ldots , A_m\}$ and associate the 
system with this family. We also assume a basis in~$\re^d$ to be fixed and associate 
matrices with the corresponding linear operators.  
By $\cA^k$ we denote the set of all products of matrices from~$\cA$
of length~$k$ (without ordering and with the repetitions permitted); 
$\cA^{\n}$ denotes the set of all products of lengths~$k\ge 1$.  
\begin{defi}\label{d.jsr}
The joint spectral radius (JSR) of a family~$\cA$ is 
 \begin{equation}\label{eq.jsr}
 \rho(\cA) \ = \ \lim_{n \to \infty}\max_{\Pi \in \cA^n} \, \|\Pi\|^{1/n}\, .
 \end{equation} 
 \end{defi}
The limit in~(\ref{eq.jsr}) always exists and does not depend on the matrix 
norm~\cite{RS}. 
For one matrix~$\cA = \{A\}$, the JSR becomes the usual spectral radius~$\rho(A)$, 
which is the largest modulus of its eigenvalues, i.e., the modulus of a leading eigenvalue. 
JSR has been studies in the literature due to numerous  
applications~(see bibliography in~\cite{J}).  

As usual, we define a {\em convex body} in~$\re^d$
as a convex compact set with a nonempty interior. If the converse is not stated, we always assume  convex bodies and polytopes to be symmetric about the origin.  
 \begin{defi}\label{d.invar}
A convex body $G\subset \re^d$ is called invariant for a matrix family~$\cA$
if it satisfies equation~(\ref{eq.invar}). 
 \end{defi}
The existence of an invariant body for any irreducible matrix family 
was proved by A.Dranishnikov and S.Konyagin in 1993 and was first published 
in 1996~\cite{P96} with a new proof. Then in~\cite{PW} it was  shown that 
a polar to an invariant body of the transpose family~$\cA^* =  \{A_1^T, \ldots , A_m^T\}$
is a unit ball of the Barabanov norm for~$\cA$. Thus, there is a one-to-one correspondence between invariant bodies and Barabanov norms. Therefore, we will formulate our results 
for both those objects. 

To a word~$s_1\ldots s_k$ of the alphabet~$\{1, \ldots , m\}$, we associate
the product~$A_{s_k}\ldots A_{s_1} \in \cA^k$.  Note that the order of multipliers is 
inverse to the order of letters!  
A prefix is some left subword of the word and a suffix is a right subword. 
The product of several words is their
concatenation. 

We use the trigonometric form of the complex number $z = |z|e^{\,\varphi i}$, 
where $\varphi$ is the argument of~$z$. If $\frac{\varphi}{\pi} \in \rat$, then we say that 
$z$ has a rational $\, {\rm mod}\, \pi\, $  argument.

For an arbitrary convex body~$\, G \subset \re^d$ symmetric about the origin, 
$\| \cdot  \|_G$ denotes the Minkowski norm~$\| \bx  \|_G \, = \, 
\sup \, \bigl\{ \lambda \ | \ \lambda^{-1} \bx \in G\, \bigr\}$.

We denote by $\Gamma$ the unit circle on the two-dimensional plane and by $D$ the unit disc.  
For an arbitrary  set $K \subset \re^2$ and arbitrary vectors 
$\bx, \by \in \re^d$, we denote by $\Phi_{\bx, \by}(K)$ 
the image of $K$ under the map $\re^2 \to \re^d$ that takes 
the basis of~$\re^2$ to vectors $\bx, \by$. This map is given by the $d\times 2$ matrix 
composed of two columns $\bx, \by$. In particular, 
$\Phi_{\bx, \by}(\Gamma)$ is an ellipse. If a matrix $A$ has a complex leading 
eigenvector~$\bv = \bx + i \by$, then the ellipse $\Phi_{\bx, \by}(\Gamma)$
will be called {\em leading} and its linear span is the {\em leading eigenspace}
(or  {\em leading plane}).

As usual, the asymptotic equivalence $\asymp$ means the existence of two 
positive constants $C_1, C_2$ such that $\, C_1\, \rho^k \, \le \, \|\bx(k)\| \, \le \, C_2\, \rho^k$. 

\newpage

\bigskip

\begin{center}
\large{\textbf{2. Fundamental theorems}}
\end{center}
\bigskip 

We consider the set of products $\cA^{\,\n}\, = \, 
\{A_{s_k}\ldots A_{s_1} \ | \ A_{s_i} \in \cA, \, k \in \n\}$ of a finite family of matrices~$\cA$. 
A product $\Pi \in \cA^{\,\n}$ is called {\em primitive} if it is not a power of a shorter product.
For a given product~$\Pi$, we denote its length (the number of multipliers) by~$|\Pi|$
and~$\nu(\Pi) = [\rho(\Pi)]^{1/|\Pi|}$.    
\begin{defi}\label{d.10}
Let $\cA$ be a finite family of matrices. A set~$\cP \subset \cA^{\n}$  is called a set of dominant products if 

1) all products from~$\cP$ are primitive and are all  different up to cyclic permutations, i.e., none of them is a cyclic permutation of another; 

2) there is a number~$q$ such that $\nu(\Pi) = q$ for all~$\Pi \in \cP$; 

3) there is $\varepsilon > 0$ such that for every $S \in  \cA^{\n}$, we have 
$\nu(S) \, \le \, (1-\varepsilon)^{1/|S|}\, q$ 
unless~$S$ is a power of some product from~$\cP$ or of one of its cyclic permutations.
\end{defi}  
\begin{remark}\label{r.3}
{\em At first, this definition  may seem impossible to verify 
within finite time since part (3) involves infinitely many conditions. 
Nevertheless, this  can be done efficiently by the {\em invariant polytope algorithm} presented in~\cite{GP1, GP2}. Given candidate products are dominant if and only 
if  this algorithm halts. 
We recall the algorithm and discuss this issue in detail in Section~4.   }
\end{remark}

If $\cP$  is a set of dominant products for the family~$\cA$, then 
$q$ is equal to the joint spectral radius~$\rho(\cA)$, see~\cite{GP1}. Moreover, 
 for all matrix products~$S \in \cA^{\,\n}$, the value $q^{-|S|}\rho(S)$ is either 
equal to~$1$ (if $S$ is a power of a dominant product or of one of its cyclic permutations), or is at most~$1-\varepsilon$. Hence, the interval 
$\Bigl( 1-\varepsilon\, , \, 1\Bigr)$ can be called a {\em spectral gap}:  
no numbers~$q^{-|S|}\rho(S), \, S\in \cA^{\,\n},$ belong to it.  

If $q=1$, then the dominance property can be defined in a simpler way: 
there is $\varepsilon > 0$ such that, for every product $S \in  \cA^{\n}$, we have 
$\rho(S) \, \le \, 1-\varepsilon\,$ 
unless~$S$ is a power of some product from~$\cP$ or of one of its cyclic permutations, 
in which case~$\rho(S) = 1$.
Hence, there is an equivalent definition of the dominant set: the set is dominant 
if 1) and 2) hold and for the normalized family 
$\tilde \cA = \{ \tilde A_i \, = \, q^{-1}A_1, \ i = 1, \ldots , m\}$,  
there is $\varepsilon > 0$ such that for every $\tilde S \in  \tilde \cA^{\n}$, we have 
$\rho(\tilde S) \, \le \, 1-\varepsilon$ 
unless the corresponding product~$S$ is a power of some product from~$\cP$ or of one of  its cyclic permutations. 
This way the dominance has been defined in~\cite{GP1} for one product and then extended for arbitrary set of products in~\cite{GP2}.  


We mostly deal with two cases. If $\cP = \{\Pi\}$ is a one-element set, we say that 
$\Pi$ is a dominant product (always assuming that it is unique).
In this case we say that~$\cA$ has a unique dominant product, although it is 
actually unique only up to a cyclic permutation.  
If $\cP = \{\Pi^{(1)}, \ldots , \Pi^{(r)}\}$ is a finite set, then we say that 
the family~$\cA$ has finitely many dominant products. 

We always make an assumption that each dominant product~$\Pi^{(i)}$ 
has a unique and simple leading eigenvalue $\lambda$. 
This means that $\lambda$ is not multiple and all other eigenvalues 
(except for the complex conjugate~$\bar \lambda$ if $\lambda\notin \re$) are strictly smaller than 
$\lambda$ in modulus. 

Let us recall that by uniqueness of an invariant body or of a norm we always 
mean their uniqueness up to multiplication by a constant. 

\begin{theorem}\label{th.10}
   Let a family of operators~$\cA$ have a unique 
   dominant product with a unique and simple  leading eigenvalue~$\lambda$.  
   If $\lambda$ is either real or complex with an irrational $\, {\rm mod}\, \pi\, $ argument, then~$\cA$ possesses a unique invariant body. 
If~$\lambda$ is real, then this invariant body is a polytope, if 
 $\lambda$ is complex with an irrational $\, {\rm mod}\, \pi\, $ argument, then this is a convex hull of several ellipses.    
\end{theorem}

At the first sight, the assumption of Theorem~\ref{th.10} is quite restrictive: 
the family~$\cA$ must have a unique dominant product whose leading 
eigenvalue is unique and simple. It turns out, however, that 
a vast majority of matrix families satisfies it. This observation was made 
first in~\cite{GP1} and then confirmed in~\cite{GP2, PG15, M} by analysing lots  of numerical
experiments with random families and with families from applications. 
Moreover, a dominant product can be efficiently found algorithmically~\cite{GP1} and 
the same algorithm constructs 
an invariant body~\cite{GZ}. We analyse this issue in Section~4.   

What can be said in the case which is not covered by Theorem~\ref{th.10}:
when the leading eigenvalue~$\lambda$ is non-real but 
possesses a rational $\, {\rm mod}\, \pi\, $ argument? 
In this case there is still an invariant body as a convex hull of ellipses, 
but it is never unique: there exist infinitely many invariant bodies of other form.  
 \begin{prop}\label{p.5}
   Suppose a family of operators~$\cA$ has a unique 
   dominant product whose leading eigenvalue~$\lambda$ is non-real and  has  
  a rational $\, {\rm mod}\, \pi\, $   argument;  then~$\cA$ has infinitely many 
    invariant bodies, one of which is a convex hull of several ellipses.    
\end{prop}
In Section~8 we classify all invariant bodies for the case of 
non-real eigenvalue with a rational  $\, {\rm mod}\, \pi\, $ $\pi$ argument.
Note that the  transpose family $\cA^*$ possesses the same property of the 
uniqueness of the dominant product (with the same leading eigenvalue). 
Applying Theorem~\ref{th.10} and Proposition~\ref{p.5} to the transpose family
and taking the polar of the invariant body, we obtain the following theorem 
that classifies Barabanov norms for generic matrix families.

\begin{theorem}\label{th.12}
   Let a family of operators~$\cA$ have a unique 
   dominant product with a unique and simple  leading eigenvalue~$\lambda$. If~$\lambda$ 
    is  real, then $\cA$ has a unique Barabanov norm.
    This norm is 
  piecewise-linear and  is given by  the formula
 \begin{equation}\label{eq.lin}
 f(\bx) \quad = \quad \max_{\bv^*} \, \Bigl|(\bv^*\, , \, \bx)\Bigr|,
 \end{equation}
  where the maximum 
 is taken over all vertices~$\bv^*$ of the invariant polytope~$G^*$ of the transpose family~$\cA^*$.  
 If~$\lambda$ is complex, then $\cA$ has a piecewise-quadratic 
 Barabanov norm given by the formula 
 \begin{equation}\label{eq.qua}
 f(\bx) \ = \ \max_{E^*}\, \max_{\bz^* \in E^*} \, \Bigl|(\bz^*, \bx)\Bigr|,
 \end{equation}
  where the maximum 
 is taken over all ellipses~$E^*$ that form the invariant body of the transpose
 family~$\cA^*$. If the argument of~$\lambda$ is irrational $\, {\rm mod}\, \pi$, 
 then this Barabanov norm  is unique. 
\end{theorem}

The Barabanov norm~(\ref{eq.qua}) can be written in a simpler form~(\ref{eq.qua1}), 
see Remark~\ref{r.5} below. 
Thus, if $\lambda \in \re$, then the unit ball of~$f$ is a polyhedron 
which is a polar to the invariant polytope of the transpose family~$\cA^*$. 
In this case $\cA$ has no other Barabanov norms. 
If $\lambda \notin \re$, then~$f$ is piecewise-quadratic; its  unit ball is the intersection of right elliptic cylinders with two-dimensional bases, those cylinders are polars to the ellipses 
forming  the invariant body 
of the transpose family~$\cA^*$.  If the argument of~$\lambda$ is irrational $\, {\rm mod}\, \pi$, then 
$\cA$ has no other Barabanov norm. If the argument is rational $\, {\rm mod}\, \pi$ but 
$\lambda$ is non-real,  
then this norm is not unique and the family~$\cA$ has infinitely many Barabanov norms. Their complete classification is 
obtained in Section~8. 

\begin{remark}\label{r.5}
{\em If an ellipse $E^*$ is defined by a pair of vectors 
$\ba,\bb \in \re^d$, i.e., $E^* = \Phi_{\ba, \bb}(\Gamma)$, 
then  $\max_{\bz^* \in E^*} \, \Bigl|(\bz^*, \bx)\Bigr| \, = \, 
\sqrt{(\ba, \bx)^2 + (\bb, \bx)^2}$, hence the formula~(\ref{eq.qua})
for the Barabanov norm~$f(\bx)$ can be written as follows:   
 \begin{equation}\label{eq.qua1}
 f(\bx) \ = \ \max_{i = 1, \ldots , N}\, \sqrt{(\ba_i, \bx)^2 + (\bb_i, \bx)^2},
 \end{equation}
where $(\ba_i, \bb_i)$ is the pair of vectors defining the $i$th ellipse  
$E^*_i = \Phi_{\ba_i, \bb_i}(\Gamma)$ in the convex hull 
$G^* = {\rm co}\ \{E^*_1, \ldots , E^*_N\}$ for the invariant body~$G^*$ of the transpose
 family~$\cA^*$. 
}
\end{remark}

Proofs to Theorems~\ref{th.10} and~\ref{th.12} 
are given in Section~5.   In Section~4 we address the practical issue: 
how to prove how to prove that the assumptions of those theorems are satisfied 
and how to construct the invariant body and the Barabanov norm. 
Now we give several illustrative examples in dimensions $d=2$ and $d=3$, 
with the corresponding pictures of invariant sets and unit balls of the Barabanov norms. 
Numerical results for higher dimensions (of course, without pictures) are considered 
later in Section~10. 

 \begin{remark}\label{r.110} Comparison with known results on the uniqueness of the Barabanov norm. 
 {\em Various sufficient conditions for the  uniqueness of  Barabanov's norm 
 have been proposed in~\cite{Ma08, M1, M2}. In~\cite{M1} it was shown that if 
 the family~$\cA$ of matrices satisfies the so-called} unbounded
agreements {\em and possesses the} rank one property, {\em then Barabanov's norm is unique.  
Both conditions  are hard to verify apart from special cases. 
However, Theorem~\ref{th.120} proved below in Section~5 implies that 
if $\cA$ possesses a unique dominant product with real and simple 
leading eigenvalue, then both those conditions are satisfied. 
So, in this special case the the main result of~\cite{M1}
implies the uniqueness part of Theorem~\ref{th.10}. Although this implication is not 
straightforward and requires a proof using Theorem~\ref{th.120}. 

Another sufficient  uniqueness condition presented in~\cite{M2}, the} transitivity property {\em seems to be very 
particular. For instance, under the assumption of Theorem~\ref{th.10} it is never satisfied 
for dimensions~$d > 2$.   

 The uniqueness issue was addressed in~\cite{Ma08} but no corresponding results have been
 obtained there.    
 }
 \end{remark}

\bigskip 

\newpage 

\begin{center}
\large{\textbf{3. Examples}}
\end{center}
\bigskip

We consider several low-dimensional examples illustrating 
Theorems~\ref{th.10} and~\ref{th.12}. 
In all the cases
the computation took a few seconds an a standard laptop.  
Higher dimensions (up to 20 for general matrices and to 2000 for nonnegative matrices)
are addressed in Section~10. 

\begin{ex}\label{ex.10}
{\em For the family $\cA = \{A_1, A_2\}$, where 
\begin{equation}\label{eq.ex10}
A_1 \ = \ 
\left(
\begin{array}{rr}
2 & -2\\
1 & 2
\end{array}
\right)\ ; \qquad 
A_2 \ = \ 
\left(
\begin{array}{rr}
1 & 2\\
-1 & -3
\end{array}
\right)\ , 
\end{equation}
the dominant product is $\Pi = A_1^3A_2$, the leading eigenvalue is real. 
The invariant convex body~$G$ 
is a $10$-gon  (Fig.~1, left). 
Its polar $G' \, = \, \{\bx \in \re^d \ | \ \max_{\by \in G^*} (\by, \bx) \le 1\}$
is the unit ball for the (unique!)  Barabanov norm for
the transpose family~$\cA^* = \{A_1^T, A_2^T\}$, Fig.~1 (right). It is also a $10$-gon.

\begin{figure}[ht!]
\begin{minipage}[h]{0.8\linewidth}
\center{\includegraphics[width=1\linewidth]{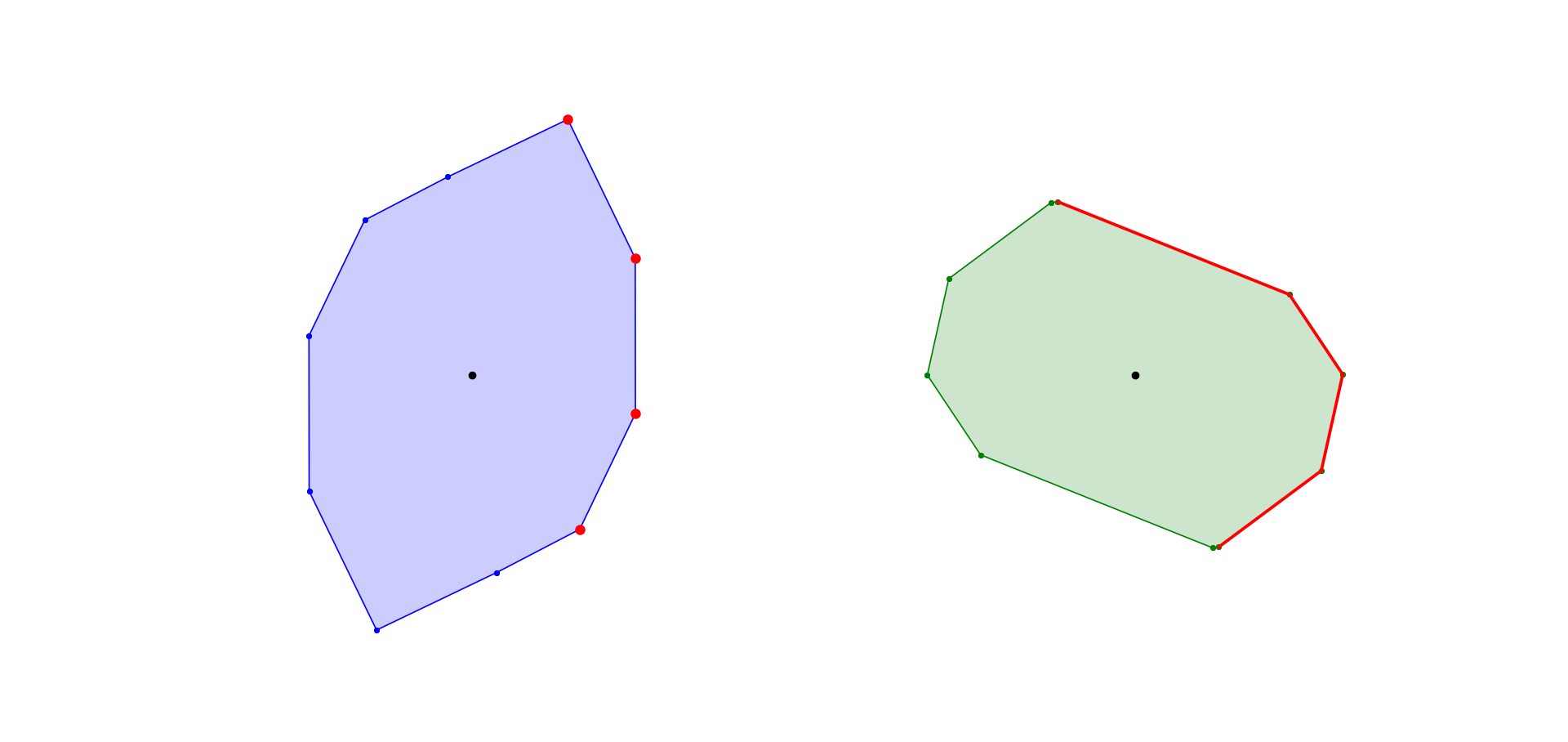}} 
 
\end{minipage}
\caption{{\footnotesize Real case, $d=2$. \textbf{Left:} the invariant polygon $G$ for 
the family~(\ref{eq.ex10}); 
\textbf{Right:} its polar $G'$ is the unit ball for the Barabanov norm of
the transpose family $\cA^* = \{A_1^T, A_2^T\}$.}}
\label{fig10}
\end{figure}
}
\end{ex}

\begin{ex}\label{ex.20}
{\em The family 
\begin{equation}\label{eq.ex20}  
A_1 \ = \ 
\left(
\begin{array}{rrr}
1 & 2& 1\\
-1 & 3& 2\\
2& -2& 3
\end{array}
\right)\ ; \qquad 
A_2 \ = \ 
\left(
\begin{array}{rrr}
-1 & 0& 3\\
0 & -1& -2\\
-3& 2&1
\end{array}
\right)\ 
\end{equation}
has a dominant product $\Pi = A_1^2A_2$ with a real leading eigenvalue. 
The invariant polytope~$G$ 
has $24$ vertices and $44$ faces (Fig.~2, left).  

Its polar $G'$ 
is the unit ball for the (unique!) Barabanov norm of~$\cA^*$, Fig.~2 (right). 
This is a polytope with $44$ vertices and $24$
faces.

\begin{figure}[ht!]
\begin{minipage}[h]{0.8\linewidth}
\center{\includegraphics[width=1\linewidth]{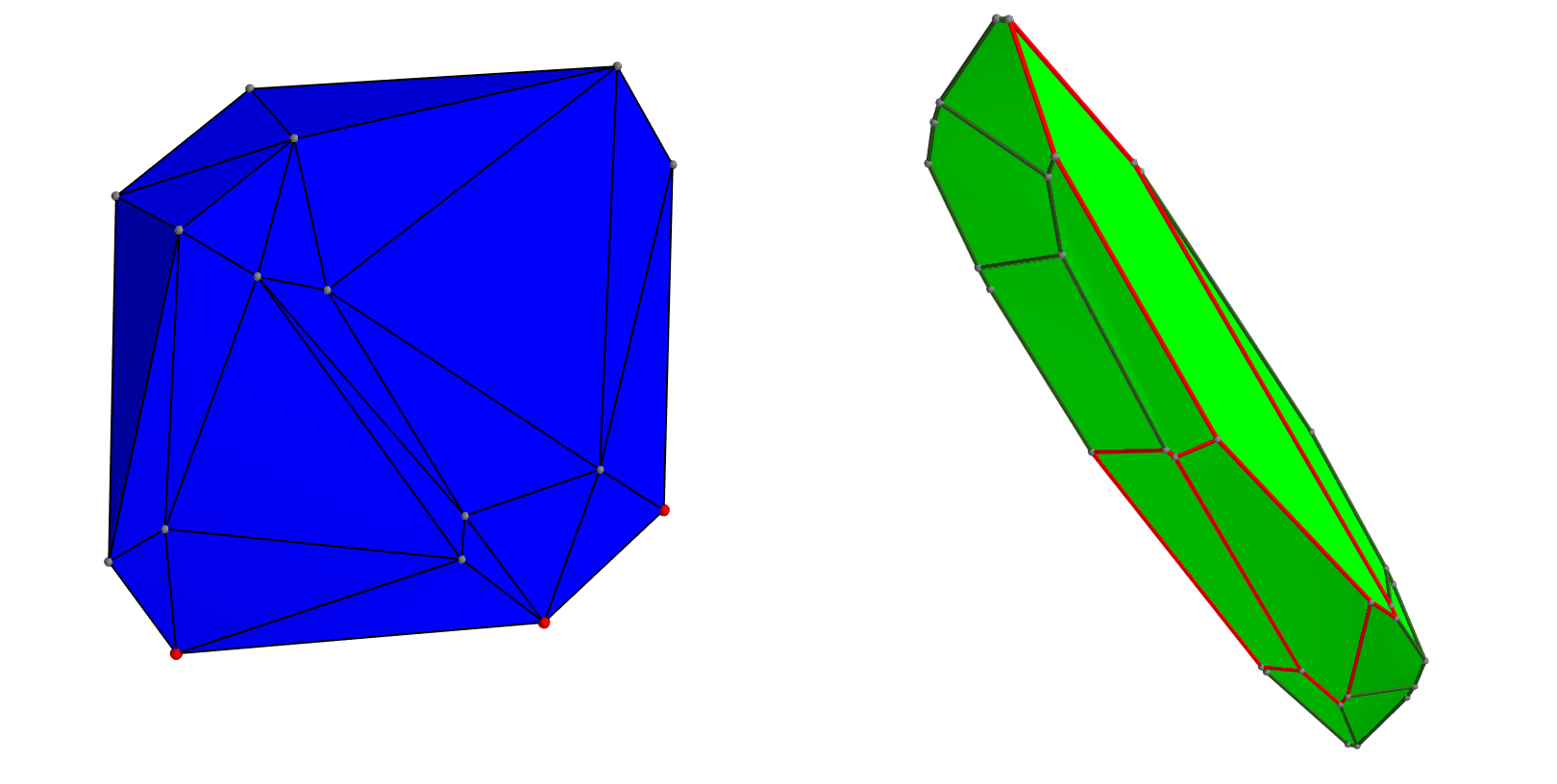}} 
 
\end{minipage}
\caption{{\footnotesize Real case, $d=3$.  \textbf{Left:} the invariant polytope $G$ for 
the family~(\ref{eq.ex20}); 
\textbf{Right:} its polar  $G'$, which is the unit ball for the Barabanov norm of $\cA$.}}
\label{fig20}
\end{figure}
}
\end{ex}

\begin{ex}\label{ex.30}
{\em The family 
\begin{equation}\label{eq.ex30} 
A_1 \ = \ 
\left(
\begin{array}{rr}
0 & 1\\
-1 & 0
\end{array}
\right)\ ; \qquad 
A_2 \ = \ 
\left(
\begin{array}{rr}
0.890 & 0.646\\
-0.129 & -0.178
\end{array}
\right)\ 
\end{equation}
has a dominant product $\Pi = A_1$ with a complex leading eigenvalue. 
The invariant convex body~$G$ 
is a convex hull of three ellipses (Fig.~3, left). 
Its polar is the intersection of three ellipses, it is the unit ball for the Barabanov norm for~$\cA^*$ (Fig.~3, right).

\begin{figure}[ht!]
\begin{minipage}[h]{0.8\linewidth}
\center{\includegraphics[width=1\linewidth]{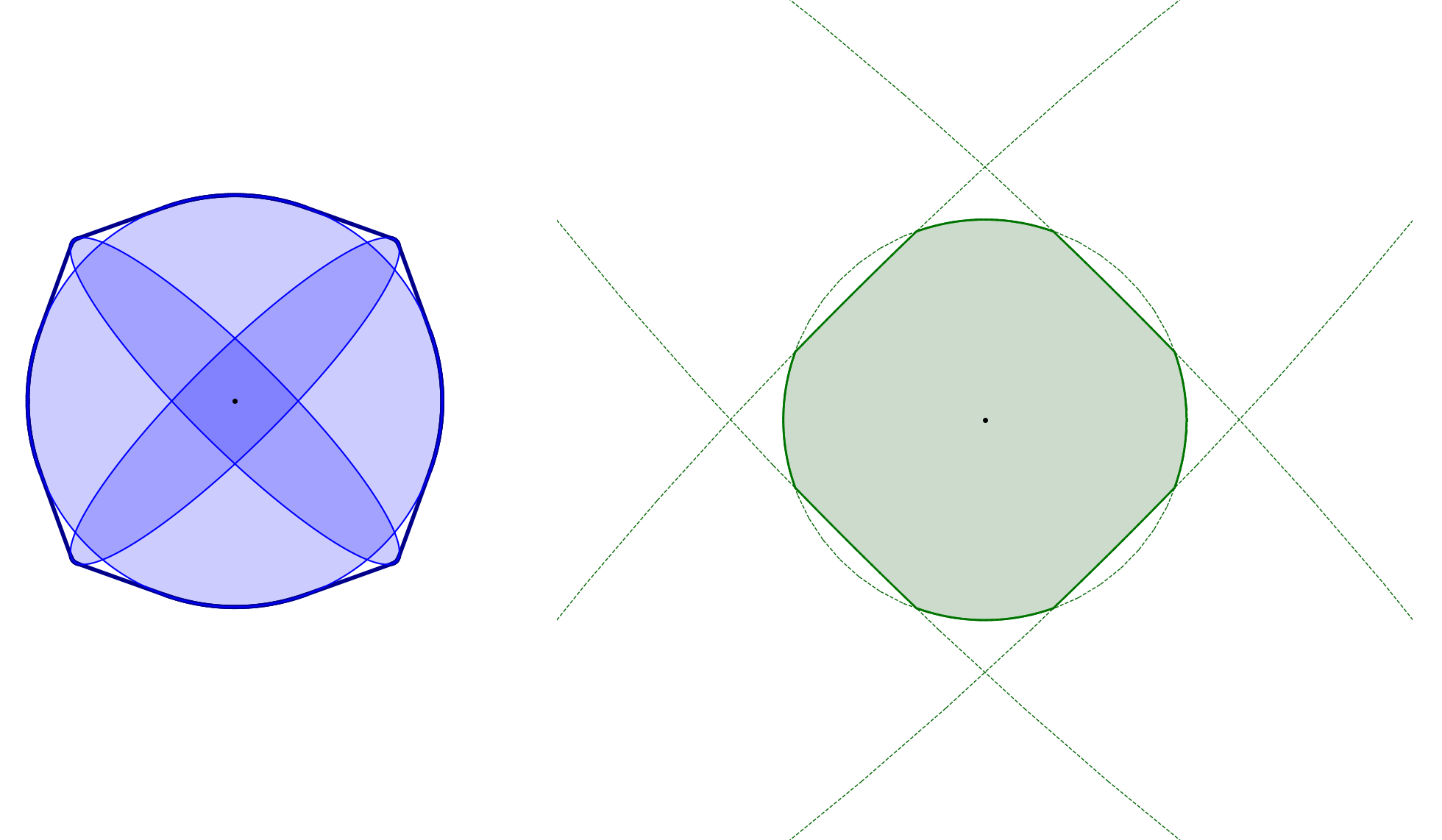}} 
 
\end{minipage}
\caption{{\footnotesize Complex case, $d=2$. \textbf{Left:} the invariant body $G$ for 
the family~(\ref{eq.ex30}); 
\textbf{Right:} its polar  $G'$.}}
\label{fig30}
\end{figure}
}
\end{ex}

\begin{ex}\label{ex.40}
{\em The family 
\begin{equation}\label{eq.ex40}
A_1 \ = \ 
\left(
\begin{array}{rr}
0 & 1\\
-1 & 0
\end{array}
\right)\ ; \qquad 
A_2 \ = \ 
\left(
\begin{array}{rr}
0.340 & 1.046\\
-0.523 & 0.170
\end{array}
\right)\ 
\end{equation}
has a  dominant product $\Pi = A_1$, with a complex  leading eigenvalue. 
The invariant convex body~$G$ 
is a convex hull of 9 ellipses (Fig.~4, left). 
Its polar is the intersection of 9 ellipses, it is the unit ball for the Barabanov norm for~$\cA^*$ (Fig.~4, right).

\begin{figure}[ht!]
\begin{minipage}[h]{0.8\linewidth}
\center{\includegraphics[width=1\linewidth]{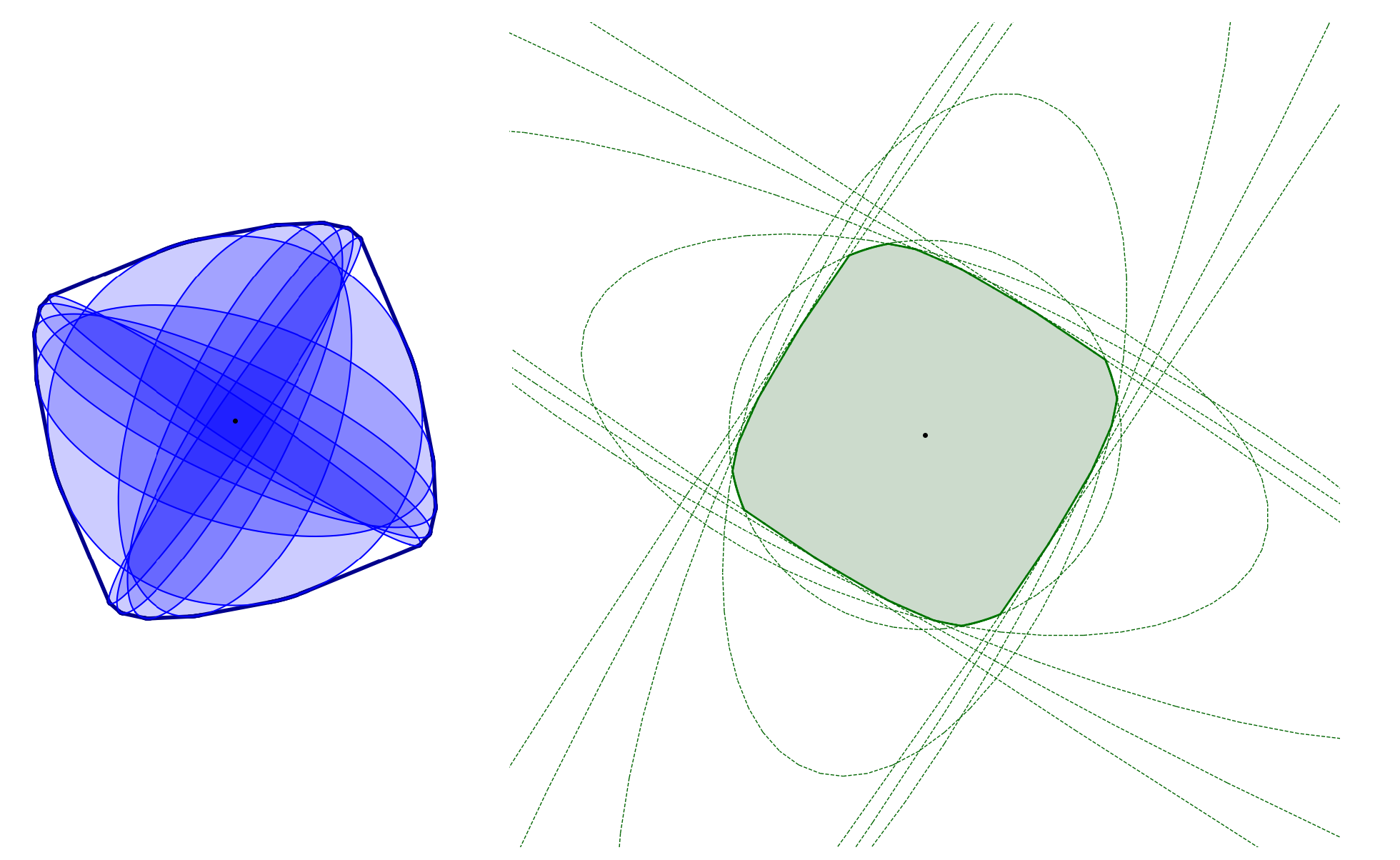}} 
 
\end{minipage}
\caption{{\footnotesize Complex case, $d=2$. \textbf{Left:} the invariant body $G$ for 
the family~(\ref{eq.ex40}); 
\textbf{Right:} its polar  $G^'$.}}
\label{fig40}
\end{figure}
}
\end{ex}

\begin{ex}\label{ex.45}
{\em The family of $3\times 3$ matrices 
\begin{equation}\label{eq.ex45}
A_1 \ = \ 
\left(
\begin{array}{rrr}
-4436 & -3993 & 887\\
3045 & -257 & -359\\
2416 & 1895  & 1338
\end{array}
\right)\ ; \ 
A_2 \ = \ 
\left(
\begin{array}{rrr}
2598 & 2948 & 682  \\
-1424& -4331 & 2691\\
  821 & -1390 & -388
\end{array}
\right)\  
\end{equation}
has a  dominant product $\Pi = A_1$, with a complex  leading eigenvalue. 
The invariant convex body~$G$ 
is a convex hull of 6 ellipses (Fig.~\ref{fig45}).

\begin{figure}[ht!]
\begin{minipage}[h]{0.7\linewidth}
\center{\includegraphics[width=1\linewidth]{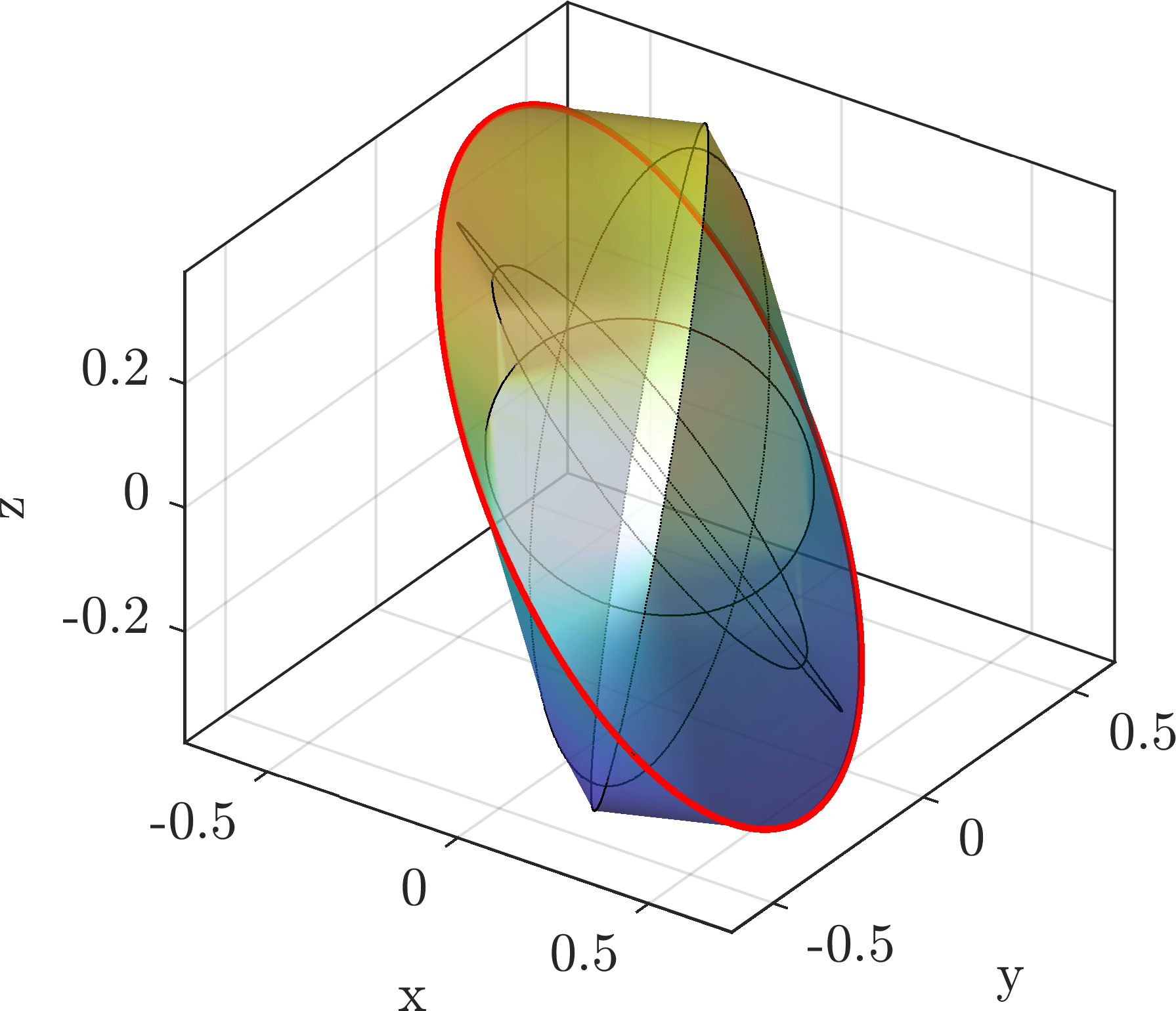}} 
\end{minipage}
\caption{{\footnotesize Complex case, $d=3$. The invariant body $G$ for 
the family~(\ref{eq.ex45})}}
\label{fig45}
\end{figure}
}
\end{ex}

\bigskip

\begin{center}
\large{\textbf{4. Construction of the invariant body and of  the Barabanov norm}}
\end{center}
\bigskip

Theorems~\ref{th.10} and \ref{th.12} in Section~2 assert that the Barabanov norm is 
unique and has a simple form, provided the system has a dominant product
with a unique and simple leading eigenvalue. 
In this section we will see that this assumption is not restrictive and 
is fulfilled for a vast majority of matrix families. Moreover, this unique Barabanov norm can be constructed in an explicit form.  This is done by the invariant polytope algorithm presented in~\cite{GP1} for the computation of JSR. 
That algorithm does a search of a matrix product~$\Pi \in \cA^{\n}$
with the biggest value of~$\nu(\Pi)$ 
among all products of some bounded lengths and then rigorously approves that 
$\rho(\cA) \le \nu(\Pi)$. The opposite equality $\rho(\cA) \ge \nu(\Pi)$ 
holds for all products~$\Pi$, this is well-known~\cite{RS}. 
This implies that~$\rho(\cA) =
 \nu(\Pi)$, and the JSR is found. 
To prove that~$\rho(\cA) \le \nu(\Pi)$  the algorithm constructs 
either a polytope~$G$ or a convex hull of several ellipses
(depending on the leading eigenvalue of~$\Pi$, which can be either real or complex)
such that ${\rm co} \bigl(\cup_{A_i \in \cA} A_iG \bigr)\, \subset \, \nu(\Pi)G$, 
which proves that $\rho(\cA) = \nu(\Pi)$. 
 Of course, there is no guarantee that the algorithm terminates within finite time.  There are examples of matrix families for which such a product~$\Pi$
does not exist~\cite{BTV, Sid1}. Nevertheless, numerical experiments and applications show 
that for a vast majority of matrix families the algorithm halts within finite time and finds the required product~$\Pi$. 
The implementation details of the algorithm were upgraded in~\cite{M, GP2}. 
Now it finds the JSR and the dominant product within a reasonable time for 
matrices of dimensions up to 20-25. The computation in higher dimensions usually takes 
too long. The version of the algorithm for non-negative matrices (see Section~9) works much faster and finds the JSR  even in dimensions of several thousands. 

Later it was observed~\cite{GZ} that the polytope~$G$ produced by the 
algorithm is nothing else but the invariant body of the family~$\cA$. Moreover, from our Theorem~\ref{th.10} (Section~2) and Theorem~\ref{th.15} below in this section, it follows that~$G$ is a unique invariant body. So, having found it once by the algorithm we can be sure that there are no others. Only if the leading eigenvalue of~$\Pi$ is non-real and has a  rational $\, {\rm mod}\, \pi\, $  argument, 
then there are infinitely many 
invariant bodies. We classify them all in~Section~8 and modify the algorithm for that case. 

Now we need to briefly recall the invariant polytope algorithm, which will be referred to as Algorithm~1. In~\cite{GP1} the cases of real and complex eigenvalues were considered separately. Here we combine them in one algorithm. 
\bigskip 

\textbf{Algorithm 1}. 
\smallskip 

\textbf{I. Choosing the candidate product}. We choose a matrix product~$\Pi = A_{s_{n}}\cdots  A_{s_1}$
({\em a candidate product}) and want to prove that $\rho(\cA) = \nu(\Pi)$. 
There are several methods to select the candidate product. 
One can just exhaust all matrix products up to  some length and 
take one 
 which attains the maximal value of~$\nu(\Pi)$. 
 There are more sophisticated methods, using branch-and-bound approach, 
 etc. see~\cite{M}.  
 
 Then we normalize our matrices as follows: $\tilde A_j = [\nu(\Pi)]^{-1}A_j$, $\tilde \cA = \{\tilde A_1, \ldots , \tilde A_m\}$ and $\tilde \Pi$ is the corresponding product of matrices from~$\tilde A$. 
 \smallskip 
 
 \textbf{II. The routine.} 
  \smallskip 
  
Let $\bv$ be the leading eigenvector of~$\tilde \Pi$. 
We define the set $V_1$ as follows. If the leading eigenvalue of~$\tilde \Pi$ is real, then 
$V_1 = \{\bv, -\bv\}$. If it is complex and, respectively,  $\bv = \bx + i \by$,  with $\bx, \by \in \re^d\setminus\{0\}$, then $V_1 = \Phi_{\bx, \by}(\Gamma)$ is an ellipse. We have $\tilde \Pi V_1 = V_1$.  
Define $V_j =  \tilde A_{s_{j-1}}\cdots  \tilde A_{s_1}V_1, \,
j= 2, \ldots , n$. The products $\tilde \Pi_j \, = \, \tilde A_{s_{j-1}}\cdots  \tilde A_{s_1} \tilde A_{s_{n}}\cdots  \tilde A_{s_j}, \, j = 2, \ldots , n$
are cyclic permutations of~$\tilde \Pi$. 
If we put formally~$\tilde  \Pi_1 = \tilde \Pi$ and $V_{n+1} = V_1$, then $\tilde \Pi_j V_j = V_j$ and $\tilde A_{s_j}V_j = V_{j+1}$ for all~$j=1, \ldots , n$.  
The set $\cR = \{V_1, \ldots , V_n\}$ is called a {\em root}.
Then we construct a sequence of finite sets $\cV_i$ and their
subsets~$\cH_i \subset \cV_i$ as follows:
 \smallskip 
 
{\tt Zero iteration}. We set $\cV_0 = \cH_0 = \cR$.
 \smallskip 
 
{\tt $k$th iteration, $k \ge 1$}. We have a finite set $\cV_{k-1}$ and its subset
$\cH_{k-1}$.
We set $\cV_{k} = \cV_{k-1},  \, \cH_k = \emptyset$ and for every $V \in \cH_{k-1}, \,  \tilde A\in \tilde \cA$, check whether $ \tilde A\, V$ is in the interior of
${\rm co} \{V \ | \ V \in \cV_{k}\}$. 
If it is, then we omit the set $ \tilde A\, V$ and take the next pair $(V , \tilde \cA) \in \cH_{k-1}\times  \tilde \cA$, otherwise we add $ \tilde A \, V$ to $\cV_{k}$ and to $\cH_k$. If $k\ge 2$, we do this for all pairs from~$\cH_{k-1}\times  \tilde \cA$. If $k = 1$ and hence $\cH_{k-1} = \cH_0 = \cR$ and $V=V_j \in \cR$, then we exclude $n$ pairs~$(V_j, \tilde A_j), \, j = 1, \ldots , n$. 

When all pairs $(V,  \tilde A)$ are exhausted, both $\cV_{k}$ and $\cH_k$ are constructed.
We define  $G_k = {\rm co} \{V \ | \ V \in \cV_{k}\}$ and  have
$$
\cV_{k}
 \, 
= \, \cV_{k-1} \cup \cH_k\, , \quad
G_{k} \, = \, {\rm co}\,  \{ \tilde A_1G_{k-1}, \ldots
  ,  \tilde A_m G_{k-1}\}\, .
$$

{\tt Termination}. The algorithm halts when $\cV_{k} = \cV_{k-1}$, i.e., $\cH_k = \emptyset$ (no new sets~$V$ are
added in the $k$th iteration). In this case
$G_{k-1} = G_k$.  Hence ${\rm co}\, \bigl\{ \cup_{j=1}^m A_jG_{k}\bigr\}\, = \, 
\nu(\Pi)\, G_k$. Therefore, $G_k$  is an invariant convex body
for~$\cA$ and  $\, \rho(\cA) \,  = \, \nu(\Pi)$.
\smallskip 

\hfill \textbf{End of the algorithm}.
\bigskip

\textbf{Implementation details}. In practice the algorithm works not with the sets 
$V_i$ but  with single points (in case or real leading eigenvalue) or with pairs of 
points (in case or a complex leading eigenvalue). 
\smallskip 

{\em The case of real leading eigenvalue of~$\Pi$.} We replace $V_1$ by 
the leading eigenvector $\bv_1$ of~$\Pi$ and then in each iteration, the set $\tilde AV_j$ is replaced by 
$\tilde A\bv_j$. Thus, all $\cV_i$ and $\cH_i$ become sets of points. Let us have in some 
step $\cV_k \, = \, \{V_i\}_{i=1}^{\ell}$. 
To decide whether a newly born set~$V = \{\pm \bv\}$ lies in the interior of $G_k \, = \, 
{\rm co}\, \{\pm \bv_i\}_{i=1}^{\ell}$ we 
solve a linear programming problem  
\begin{equation}\label{eq.lp}
\left\{
\begin{array}{l}
t_0\ \to \ \max\\
\mbox{subject to:}\\
-s_i \, \le \, t_i \, \le \, s_i , \ i = 1, \ldots , \ell\\
\sum_{i=1}^{\ell} s_i \, \le \, 1\\
t_0\bv \, = \, \sum_{i=1}^{\ell} t_i \bv_i
\end{array}
\right. 
\end{equation}
We have $\bv \in {\rm int}\, G_k$ if and only if $t_0 > 1$. In practice 
we fix a small tolerance parameter 
$\delta> 0$ (usually, $\delta$ is between $10^{-8}$ and $10^{-6}$)
and decide that the set  $V$ is redundant if  $t_0 > 1+\delta$, otherwise we keep 
the points $\pm \bv$ among the vertices of an invariant polytope, although they may actually not be vertices. 
\smallskip 

{\em The case of complex leading eigenvalue of~$\Pi$.} We replace the ellipse $V_1$ by 
the pair of points $\bx_1, \by_1 \in \re^d$
such that $\bv_1 = \bx_1 + i\by_1 $ is the leading eigenvector of~$\Pi$. 
As we know, $V_1 = \Phi_{\bx_1, \by_1}(\Gamma)$. 
Then in each iteration, the ellipse $\tilde AV_j$ is replaced by the pair of points 
$(\tilde A\bx_j, \tilde A\by_j)$. Clearly, 
$V_j = \Phi_{\tilde A\bx_j, \tilde A\by_j}(\Gamma)$. 
Thus, all  $\cV_i$ and $\cH_i$ become sets of pairs of points. 
Let $(\bx, \by)$  be a newly born pair.   
To prove that  the ellipse $V= \Phi_{\bx , \by}(\Gamma)$ 
is contained in the interior of $G_k \, = \, 
{\rm co}\, \{V_i\}_{i=1}^{\ell}$ we 
solve the following  optimization problem  
\begin{equation}\label{eq.cp}
\left\{
\begin{array}{l}
t_0\ \to \ \max\\
\mbox{subject to:}\\
\sqrt{t_{k}^2 + u_k^2} \le s_k, \ k = 1, \ldots , 2\ell\\
\sum_{k=1}^{2\ell} s_k \, \le \, 1\\
t_0\bx \, = \, \sum_{i=1}^{\ell} \bigl(t_{2i-1} \bx_i \, - \, u_{2i-1} \by_i\bigr)\, + \, 
\bigl(t_{2i} \bx_i \, + \, u_{2i} \by_i\bigr)\\
t_0\by \, = \, \sum_{i=1}^{\ell} \bigl(u_{2i-1} \bx_i \, + \, t_{2i-1} \by_i\bigr)\, + \, 
\bigl(u_{2i} \bx_i \, - \, t_{2i} \by_i\bigr)
\end{array}
\right. 
\end{equation}
This is a conic programming problem and is solved by the interior point method 
on Lorentz cones (see {\tt www.mosek.com} for the corresponding software). 
If $t_0 > 1$, then 
 $V \subset {\rm int}\, G_k$. So, we remove the pair $(\bx, \by)$ if 
  $t_0 > 1+\delta$. Otherwise we set $V_{\ell+1}  = V$ 
and add this 
ellipse (i.e., the pair $(\bx_{\ell+1}, \by_{\ell+1}) = (\bx, \by)$) to both 
$\cV_k$ and $\cH_k$. Note that, in contrast to the real case, here the condition 
$t_0 > 1$ is only sufficient but not necessary for the inclusion $V \subset {\rm int}\, G_k$. 
That is why, usually the resulting set~$G$ contains many redundant ellipses~$V_i$ which are 
not ``vertices'' of~$G$, i.e., are inside $G$ and could be removed.  It slows down the algorithm but 
not significantly, see Section~10 for numerical results.

\bigskip 

\textbf{Comments and analysis of convergence}. 
Actually, the algorithm  works with the sets $\cV_k$  only, the convex bodies~$G_k$ are needed to illustrate the geometric idea. Thus, in each iteration we construct a body $G_k \subset \re^d$, which is either a polytope (in case of real eigenvalue of~$\Pi$) or a convex hull of ellipses (the case of complex 
eigenvalue), store all its vertices (ellipses)~$V_i$ 
in the set~$\cV_k$ and spot the set $\cH_k \subset \cV_k$ of newly appeared (after the previous iteration) sets~$V_i$.
Every time we check whether $ \tilde \cA G_k \subset G_k$. 
If $ \tilde \cA G_k \subset G_k$, then $G_k$ is an invariant body,
$\|\tilde A_i\|_{G_k} \le 1$ for all $i$. Otherwise, we update the sets $\cV_k$ and $\cH_k$ and continue. 
\medskip 

If Algorithm~1 halts within finite time, then the candidate product~$\Pi$
not only gives the  precise value of JSR but also is a dominant product.  

\smallskip

\noindent \textbf{Theorem A} \cite{GP1}. {\it Algorithm~1 
applied to a candidate product~$\Pi$ terminates within finite time if and only 
if~$\Pi$ is a unique dominant product for~$\cA$ and its leading eigenvalue is unique and simple.}  
\smallskip 

Thus, we can always check whether a given product is dominant or not. 
If it is, then the  invariant body of the family~$\cA$ is readily available 
as the body~$G_k$ obtained by the end of the algorithm. 
\begin{ex}\label{ex.50}
{\em For the matrices~(\ref{eq.ex10}) from Example~\ref{ex.10}, 
the dominant product is $A_1^3A_2$. The root $\cR = \{\bv_1, \bv_2, \bv_3, \bv_4\}$
consists of four vertices of the polygon $G$, they are marked in red (Fig.~1, left). 
The four corresponding  
sides of the polar~$G'$ are also red (Fig.~1, right). 

For the matrices~(\ref{eq.ex20})  from Example~\ref{ex.20}, 
the dominant product is $A_1^2A_2$. The root $\cR = \{\bv_1, \bv_2, \bv_3\}$
consists of three vertices of the polytope $G$, they are marked in red (Fig.~2, left). 
The three corresponding  
faces of the polar~$(G^*)'$ are also red (Fig.~1, right).  

For the matrices~(\ref{eq.ex30})  from Example~\ref{ex.30}, the dominant 
product is~$A_1$, the root $\cR = \{V_1\}$, where $V_1$ is a circle 
(Fig.~3, left). The same is for the matrices~(\ref{eq.ex40})  from Example~\ref{ex.40}. 
}
\end{ex}

\begin{theorem}\label{th.15}
  If  Algorithm~1 terminates after~$k$th iteration, 
  then it produces an invariant convex body~$G_k$ for the family~$\cA$. If 
  the leading eigenvalue of the product~$\Pi$
  is either real or complex with an irrational $\, {\rm mod}\, \pi\, $  argument, 
  then~$G_k$ is a unique invariant body for~$\cA$. 
\end{theorem}

The proof is given in the next section. 
Passing to the transpose family of operators~$\cA^*$
 we obtain the method of construction of Barabanov's norm 
  presented  in~\cite{GZ}. Now we can claim that there are no other
 Barabanov norms. Applying Theorem~\ref{th.15} to the family~$\cA^*$
 we obtain the following

\begin{theorem}\label{th.18}
  If  Algorithm~1 applied to the transpose family~$\cA^*$
  terminates after~$k$th iteration, 
  then it produces Barabanov's norm for the family~$\cA$. 
  
  If the leading eigenvalue~$\lambda$ of the candidate product~$\Pi^*$
  is real, then this norm is piecewise-linear, 
  and is given by formula~(\ref{eq.lin}), where $\bv^*$ runs over the set of vertices of~$G_k^*$. This is a unique Barabanov norm for~$\cA$. 
  
If $\lambda \notin \re$, then this norm is piecewise-quadratic, 
  and is given by formula~(\ref{eq.qua1}), where 
we set $\ba_i = \bx_i,  \,  \bb_i = \by_i$ and $E^*_i = \Phi_{\bx_i, \by_i}(\Gamma), \, 
i = 1, \ldots , N$, 
  runs over the set of ellipses generating~$G_k^*$. If the argument of~$\lambda$ is 
  irrational $\, {\rm mod}\, \pi$, then 
  this is a unique Barabanov norm for~$\cA$. 
\end{theorem}

The remaining case, when $\lambda \notin \re$ and the argument of~$\lambda$
is rational $\, {\rm mod}\, \pi, $ is considered in Section~8. Now we turn to the proofs of the main results.   

\vspace{1cm}

\begin{center}
\large{\textbf{5. Proofs of the fundamental theorems}}
\end{center}
\bigskip

To prove Theorems~\ref{th.10} and~\ref{th.15}  
we need one auxiliary statement on the structure of trajectories 
of an irreducible  system (Theorem~B below). 
We begin with the following well-known fact. Its proof is given for convenience of the reader. 

\begin{lemma}\label{l.10}
There is a continuous function~$\psi(\delta, z)$ on~$\re_+^2$ such that 
$\psi(0, z) = 0$ for all~$z$ and for every $d\times d$  matrix~$A$, 
the following is true: if there is a vector $\bx$ 
such that $\|A\,\bx - \bx\| \, \le \, \delta \, \|\bx\|$, 
then~$A$ has an eigenvalue~$\lambda \in \co$ such that $|\lambda - 1| \, \le \, 
\psi (\delta, \|A\|)$.   
\end{lemma}
{\tt Proof.} Without loss of generality it can be assumed that $\|\bx\| = 1$
and that $\bx = \be_1$ is the first basis vector. The polynomial $p(\lambda) \, = \, {\rm det}\, (\lambda \, I \, -\, A)$ has the leading coefficient one and 
other coefficients at most $2^d\|A\|^d$ in modulus. Moreover, 
since $\|A\,\be_1 - \be_1\| \, \le \, \delta \, $, it follows that  
the first column of the matrix $\lambda \, I \, -\, A$ has all components at most 
$\delta$ in modulus. Since the moduli of all other entries of this matrix 
 are at most $\|A\|+1$, we have $|p(1)| \, = \, |{\rm det}\, (I \, -\, A)|\, \le  \, C \, \delta$, 
 where $C \, \le \,  
 2^d (\|A\|+1)^{d-1}$. Therefore, there exists a root of~$p$ 
 on the distance at most $C_0(C \, \delta)^{1/d}$  from the number~$1$.

{\hfill $\Box$}
\smallskip

Now we turn to the structure of trajectories. Let~$\cA = \{A_1, \ldots , A_m\}$
be an arbitrary irreducible system. The irreducibility 
implies that $\rho(\cA) > 0$ and that $\cA$ possesses at least one  invariant body~$G$~\cite{P96}. After normalization it can be assumed that $\rho(\cA) = 1$,
all invariant bodies stay the same.  
A point $\bx \in \re^d$ is called~{\em recurrent} if it belongs 
to the  boundary of the invariant body~$G$ and there is a trajectory 
$\{\bx_k\}_{k\ge 0}$ such that $\bx_0 = \bx$ and some 
subsequence~$\{\bx_{k_j}\}_{j \in \n}$ tends to $\bx$ as $j \to \infty$.  

An {\em orbit} of a point~$\bx$ is the set~$\{\Pi \bx \ | \ \Pi \in \cA^{\n}\}$, 
i.e., is a union of all trajectories starting at~$\bx$. Observe that if~$\bx$ is recurrent, then the points of the trajectory of~$\bx$ are not necessarily recurrent. 
\smallskip 

\noindent \textbf{Theorem B} \cite{P96}. {\it Let $\cA$ be an irreducible family with 
$\rho(\cA) = 1$. Then for every invariant body~$G$ of~$\cA$, there exists a compact 
subset~$\bC$ of the set of recurrent points such that~$G$ is the closed convex hull of orbits 
of points from~$\bC$.} 
\smallskip 

Note that the set $\bC$ in Theorem~B depends on~$G$. 
That is why Theorem~B does not imply the uniqueness of the invariant body. 
In fact, Theorem~B holds also in cases when the invariant body is not unique, 
for instance, when the leading eigenvalue is non-real and has a 
rational $\, {\rm mod}\, \pi\, $  argument.

If a matrix $A$ has a unique simple leading eigenvalue, 
then it has a {\em leading eigenspace} which is either one-dimensional 
(the linear span of the real leading eigenvector) or two-dimensional 
(the real linear span of the real and complex part 
of the  leading eigenvector). 
\begin{prop}\label{p.10}
If a family of operators has finitely many dominant products and each of them 
 has a unique and simple leading eigenvalue, then every recurrent point 
of this family belongs to the leading eigenspace of one of these products  or of one of its cyclic permutations. 
\end{prop}
{\tt Proof.} Let the spectral gap be the interval $(1-\varepsilon, 1)$, where 
$\varepsilon > 0$.  Since the family $\cA$ is irreducible, the norms of 
all products of its matrices  are bounded by some constant~$C$.  Choose small $\delta > 0$ so that 
$\psi (\delta , C) < \varepsilon$, where the function~$\psi$ is defined in Lemma~\ref{l.10}. 
If $\bx\ne 0$ is a recurrent point, then 
there is a product $S \in \cA^k$ such that 
$\|S \bx - \bx\| \, \le \,  \delta \, \|\bx\|$. 
By Lemma~\ref{l.10}, this implies that $S$
has an eigenvalue~$\mu$ such that $\, |1-\mu| \, \le \,   
\psi(\delta , C) \, < \, \varepsilon$. Hence, 
$\mu \in (1-\varepsilon,1+\varepsilon)$. On the other hand, because of the spectral gap, 
 $\mu$ can neither be in~$(1-\varepsilon,1)$ nor bigger than one. 
 Consequently, $\mu = 1$ and therefore   $S$
is a power of some dominant product or of one of its cyclic permutations.  
It can be assumed that this is a power of a dominant product, the case of a cyclic permutation is literally the same. 
 Thus, for the point~$\bx$, there is a dominant product~$\Pi$ 
 and a sequence of integers $\{j_k\}_{k\in \n}$ such that 
  $\|\Pi^{j_k} \bx - \bx\| \to 0$ as $k \to \infty$.   
Since the leading eigenvalue of~$\Pi$ is simple, it follows that 
$\Pi^{j_k}\bx$ converges to the projection of~$\bx$ to the line containing the 
leading eigenvector (to the leading eigenspace in the complex case). 
Therefore, $\bx$ coincides with this projection, and so 
$\bx$ is the leading eigenvector (respectively, belongs to the leading eigenspace).

{\hfill $\Box$}
\smallskip

\bigskip 

{\tt Proof of Theorem~\ref{th.10}.} 
Let $G$ be an arbitrary invariant body for~$\cA$, 
$\Pi =  A_{j_n}\cdots A_{j_1}$ be the dominant product, 
$\lambda$ be its leading eigenvalue, 
 $\Pi_{k}$ be its $k$th cyclic permutation. After normalization it can be assumed 
 that $\rho(\Pi) = 1$. 

 \textbf{The case $\lambda \in \re$}. In this case $\lambda = \pm 1$ 
and we assume $\lambda = 1$, the other case is considered in the same way.  
 Denote  by~$\bv_1$ the leading 
eigenvector of $\Pi$ that belongs to~$\partial G$ (any of the two vectors). 
Then~$\bv_k = A_{j_{k-1}}\cdots A_{j_1}\bv_{1}\, , \ k = 2, \ldots , n$. 
Clearly, $\bv_k$ is the leading eigenvector of~$\Pi_k$. Moreover, $\bv_k \in \partial G$ for all $k$. Indeed,  $\|\bv_k\|_{G} =  \|A_{j_{k-1}}\cdots A_{j_1}\bv_{1}\|_{G} \le 
\|\bv_1\|_{G} \, = \, 1$, 
because the norm $\|\cdot \|_{G}$, as a Lyapunov function of the system, is non-increasing on any trajectory;  
$\|\bv_1\|_{G} = 1$ because $\bv_1 \in \partial G$. Thus, $\|\bv_k\|_{G} \le 1$. 
On the other hand, 
denoting $\bv_{n+1} = \bv_1$, we get 
$1 = \|\bx_{n+1}\|_{G} =  \|A_{j_n}\cdots A_{j_k}\bv_k\|_{G} \le 
\|\bv_k\|_{G}$ and hence $\|\bv_k\|_{G} \ge 1$. Thus 
$\|\bv_k\|_{G} = 1$ and consequently $\bv_k \in \partial G$
for all~$k$. Theorem~B and Proposition~\ref{p.10} imply that 
$G$ is the closure of the convex hull of all trajectories starting at the points~$\bv_k$. 
 Hence it is obtained by Algorithm~1 from the candidate product~$\Pi$.

\textbf{The case $\lambda \notin \re$}. The leading ellipse is $E\, =\, \Phi_{\bx, \by}\, (\Gamma)$, 
where $\bx + i\by$ is the complex leading eigenvector of~$\Pi$. 
Normalize this vector so that
 $\bx \in \partial G$.   If the argument $\varphi$ of~$\lambda$ is 
 irrational $\, {\rm mod}\, \pi,$  
then the set $\{\Pi^{k}\bx\}_{k\in \n}$ is everywhere dense on 
 $E$. 
 The sequence of norms $\{\|\Pi^{k}\bx\|_{G}\}_{k\in \n}$
 is non-increasing and it has the number $\|\bx\|_{G} = 1$ as a limit point, hence 
it is an identical one. Therefore, the ellipse $E$ lies on the boundary of~$G$
and $E$ is an intersection of this surface with the leading eigenspace 
of~$\Pi$.     
Combining Theorem~B and Proposition~\ref{p.10} 
and taking into account that all points from~$E$ are recurrent
(because for every~$\bz \in E$, the sequence 
$\{\Pi^k\bz\}_{k \in \n}$ has a limit point~$\bz$), we conclude that 
$G$ is the closure of convex hulls of all trajectories starting at  points from~$E$. 
Therefore, $G$ is a convex hull of images of~$E$ under the action of all 
products of operators from~$\cA$. Hence it is obtained  by Algorithm~1
from the candidate  product~$\Pi$.

{\hfill $\Box$}
\smallskip

{\tt Proof of Theorem~\ref{th.15}}. After the $k$th iteration we obtain the body~$G_k$, 
which is a convex hull of sets~$V \in \cV_k$, 
where $\cV_k = \cH_{0}\cup \cH_{1} \cup \cdots \cup  \cH_{k}$.
Consider the multivalued operator~$\overline{A}$ which maps every 
element~$V \in \cV_k$ to the set of elements $\tilde A_1V, \ldots , \tilde A_mV$. 
Then $\overline{A}$ maps the root $\cH_0 = \cR$ to the 
union~$\{\cH_0, \cH_1\}$ and each set~$\cH_j, \, j\ge 2,$ to $\, \cH_{j+1}$. 
Therefore, $\overline{A} \cV_{k-1} = \cV_{k}$. However, the algorithm terminates
after the $k$th iteration, hence $\cV_{k-1} = \cV_k$. 
Thus, $\overline{A} \cV_{k} = \cV_{k}$ and therefore 
${\rm co}\, \bigl\{\tilde A_1G_k, \ldots , \tilde A_mG_k\bigr\} \,  = \, G_k$. So, $G_k$ is an invariant body. 
By Theorem~A, the product $\Pi$ is dominant. Hence, we can use 
Theorem~\ref{th.10}, which implies the uniqueness of the invariant body 
in cases of real leading eigenvalue and of complex leading eigenvalue with
an irrational $\, {\rm mod}\, \pi\, $  argument.

{\hfill $\Box$}
\smallskip

Thus, we have proved the uniqueness and have established the structure of Barabanov's norms 
for general matrix families possessing dominant products with the leading eigenvalue which is either real or  complex with an irrational $\, {\rm mod}\, \pi\, $ argument. 
 The remaining case when the leading  eigenvalue is non-real and has a 
rational $\, {\rm mod}\, \pi\, $ argument is more delicate; we attack it in Section~8. To this end we need an auxiliary result characterising the growth of trajectories of an arbitrary system with $\rho(\cA) = 1$. This result is, probably, of some 
independent interest and we put in a separate section (Section~7). 

Now we are going to 
analyse systems with several different (up to powers and cyclic permutations) dominant products. This case is rather special but it plays an important role in some applications.    

\bigskip 

\begin{center}
\large{\textbf{6. Systems with finitely many dominant products}}
\end{center}
\bigskip

According to numerical experiments, the uniqueness of the 
 dominant product takes place  for almost all matrix families
(at least, randomly generated ones).
Nevertheless, in applications it happens that there are several
dominant products. It occurs  when there are some relations between matrices of the 
family. For example, in the computation  of the H\"older regularity of wavelets
and of limit functions of subdivision schemes, one needs to 
find the JSR of two special matrices $T_0, T_1$, which are sometimes both dominant, see~\cite{G, CHM, GP2} and references therein. 
A similar situation occurs in some  problems  of combinatorics, number theory, and 
formal languages~\cite{BCJ, JPB, M, MOS, P17}.    

In fact, in the results of Section~2-4 the uniqueness of the dominant product is not a restriction. Explicit classification and construction of the invariant body and of the Barabanov norms
can be realised in a similar way when the system has finitely many dominant products.
The only difference is that, as we are going to see, the Barabanov norm {\em 
is never unique in this case}: any system with several dominant products has infinitely many 
invariant bodies and Barabanov norms, which can, nevertheless, be classified (Corollary~\ref{c.10}). 
 
 The algorithm of computing the JSR for families with several 
 dominant products was elaborated in~\cite{GP2}. It is very similar 
 to Algorithm~1, but it starts with several roots~$\cR^{(j)}$
 (each root is associated to the corresponding dominant product~$\Pi^{(j)}$). 
 However, to provide the convergence of the algorithm one needs to {\em balance} 
 the roots, i.e., to multiply each of them by a certain positive constant~$\alpha_j$, and those constants have to be found.  
 Otherwise, the algorithm does not converge within finite time.  

 \bigskip

\textbf{Algorithm 2}. 

\textbf{Choosing the candidate products}. We choose 
several candidate products $\Pi^{(1)}, \ldots , \Pi^{(r)}$
that are all primitive, 
different up to cyclic permutations, and having the same values~$\nu(\Pi^{(i)})$
for all~$i = 1, \ldots , r$.   Denote this value  by $q$. 
 Then normalise all matrices from~$\cA$ as~$\tilde A_j = q^{-1}A_j$. 
 \smallskip 
 
 \textbf{The balancing}. 
Take a vector of positive numbers $\balpha = (\alpha_1, \ldots , \alpha_r)$
called  {\em balancing vector}. Those numbers are selected in a special way to provide the convergence of the algorithm. For a method of finding a proper balancing vector see~\cite{GP2}. Then we define the sets $V_1^{(i)}, \ldots ,  V_{n_i}^{(i)}$
 as in Algorithm~1: they are either symmetric pairs of leading eigenvectors
 of $\Pi^{(i)}$ (if~$\Pi^{(i)}$ has a real leading 
 eigenvalue) or of leading ellipses (the case of complex leading eigenvalue). 
 Then 
 we form the roots $\cR^{(i)} = \{\alpha_i V_1^{(i)}, \ldots , \alpha_i V_{n_i}^{(i)}\}$ (each element $V_k^{(i)}$ is multiplied by~$\alpha_i$)
 that consist of symmetric pairs of leading eigenvectors (the case of real leading 
 eigenvalue) or of leading ellipses (non-real leading eigenvalue). 
    \smallskip 
 
 \textbf{The routine.} 
 \smallskip 
 
{\tt Zero iteration}. We set $\cV_0 = \cH_0 = \cup_{i=1}^r \alpha_i\cR^{(i)}$. 
\smallskip 

{\tt $k$th iteration} is literally the same as in Algorithm~1. 
\smallskip 

{\tt Termination} is the same as in Algorithm~1. 
\smallskip 

\hfill \textbf{End of the algorithm}.
\smallskip

If Algorithm~2 terminates after the~$k$th iteration, then 
 we obtain a convex body~$G_k = {\rm co}\, \{V \ | \ 
V \in \cV_k\}$ which is an invariant body for~$\cA$. It is a convex hull of several
segments and of several ellipses, all centered at the origin.  
\medskip

\noindent \textbf{Theorem C} \cite{GP2}. {\it If Algorithm~2 
applied to  candidate products~$\Pi^{(1)}, \ldots , \Pi^{(r)}$ 
with equal values of~$\nu(\Pi^{(i)})$ 
and to some balancing vector~$\balpha$ terminates within finite time, then these products are dominant and each of them has a unique and simple leading eigenvalue. Conversely, if these products are dominant for~$\cA$ and their leading eigenvalues are unique and simple, then 
 there is a balancing vector for which Algorithm~2 terminates 
within finite time}.   
\smallskip 

Thus, if $\cA$ has finitely many dominant products, then they can be found 
by Algorithm~2 along with the weights~$\{\alpha_i\}_{i=1}^{r}$. 
Note that the set of dominant products is unique but the set of weights is not.  Let us now show that 
the same algorithm gives the invariant body and, if the leading eigenvalues of all 
the dominant products are either real or complex with irrational $\, {\rm mod}\, \pi\, $ 
arguments, then all invariant bodies are exhausted by those found with Algorithm~2.

\begin{theorem}\label{th.25}
  If  Algorithm~2 halts  within finite time making~$k$ iterations,  
  then it produces an invariant convex body~$G_k$ for the family~$\cA$. 
  
  If the leading eigenvalues of all the dominant  products  of~$\cA$ are either 
  real or complex with irrational $\, {\rm mod}\, \pi\, $  arguments, then 
  every invariant body of~$\cA$ is obtained by Algorithm~2 with some 
  balancing vector. Different balancing vectors produce different  invariant bodies. 
\end{theorem}

{\tt Proof}. If Algorithm~2 halts after~$k$th iteration, then 
$G_k$ is an invariant body. Moreover, in this case 
all the products~$\Pi^{(i)}, \, i = 1, \ldots , r$ are dominant
(Theorem~C). 
 If all their leading eigenvalues  are either 
  real or complex with an irrational $\, {\rm mod}\, \pi\, $ argument, then 
  every invariant body of~$\cA$  is obtained by Algorithm~2 with some 
  balancing vector. This is proved in the same way as Theorem~\ref{th.10} by 
applying Proposition~\ref{p.10}  and Theorem~B from Section~5. Also in the same way we show that 
all the sets from the roots~$\cR^{(i)}$ lie on the boundary of the invariant body. 
Hence, changing the multipliers~$\alpha_i$ we obtain different invariant sets.

{\hfill $\Box$}
\smallskip

The key difference with the case of one dominant product is that now 
choosing different weights we get different invariant bodies. 
This is the reason of non-uniqueness of the invariant body for families 
with many dominant products. 

\begin{cor}\label{c.10}
If a family of operators has $r\ge 2$ dominant products (up to cyclic permutations), then it has infinitely many invariant bodies.  If, in addition,  the leading eigenvalues of all those dominant  products are either 
  real or complex with irrational $\, {\rm mod}\, \pi\, $ arguments, then all those invariant bodies are convex hulls 
  of finitely many points and ellipses. 
\end{cor}
{\tt Proof}. By Theorem~C, there exists a balancing vector~$(\alpha_1, \ldots, \alpha_r)$ for which Algorithm~2 terminates within finite time and gives an
invariant body. 
Fix $\alpha_1$. If we slightly vary other coefficients~$\alpha_2, \ldots , \alpha_r$, then Algorithm~2 performs the same iterations as before. Indeed, each iterations is defined by 
the set of dead vertices, when $\tilde A V_p \in {\rm int}\, G_k$. 
A sufficiently small variation of parameters keeps this inclusion. 
Hence, after a small variation of~$\alpha_2, \ldots , \alpha_r$, Algorithm~2 
performs the same iterations. Consequently it  terminates within finite time.  By Theorem~\ref{th.25}, different variations of parameters produce different  
invariant body.

{\hfill $\Box$}
\smallskip

\begin{remark}\label{r.30}{\em 
Applying Algorithm~2 to the transpose family~$\cA^*$
we obtain the Barabanov norm for~$\cA$. Similarly to Theorem~\ref{th.18}
in Section~4
one expresses the relation between the invariant body of~$\cA^*$
and the Barabanov norm for~$\cA$.  

}
\end{remark}

Thus, in case of several dominant products the Barabanov norm is never unique. 
Nevertheless, if the leading eigenvalue of every dominant product is either real or 
complex with an irrational $\, {\rm mod}\, \pi\, $ argument, then all those Barabanov norms are classified by 
Theorem~\ref{th.25}. They are parametrized by 
the balancing vectors~$\balpha \in \re^r_+$ for which Algorithm~2 halts within finite time. Since each balancing vector can be normalized by the condition~$\sum_{i=1}^r \alpha_i = 1$,  we see that there exists a $(r-1)$-parametric family of Barabanov norms. 

If  at least one of the dominant products has a complex leading eigenvalue with a 
rational $\, {\rm mod}\, \pi\, $ argument, the Barabanov norm can still be computed by the same Algorithm~2. 
However, there will be other norms that are not obtained by that algorithm. 
Their classification requires another method, see Section~8. To introduce that method we first need to make a detailed analysis of growth of trajectories of a discrete-time system. This is a subject of the next section.

 \newpage

\bigskip 

\begin{center}
\large{\textbf{7. Classification of trajectories of the fastest growth}}
\end{center}
\bigskip

The results of this section will be applied in characterising  Barabanov norms 
in  case of rational $\, {\rm mod}\, \pi\, $ arguments of the leading eigenvalue. They are also of an independent interest. We are going to find all trajectories of a linear switching system
with the fastest asymptotic growth.

An analysis of asymptotic growth of trajectories is a 
subject of an extensive literature, see, for example,\cite{CMS, Sid2, L, M3, ZPC}. 
If the discrete-time linear switching system~(\ref{eq.sys})
is irreducible,  then the maximal possible  growth of trajectories is~$
\|\bx(k)\| \, \asymp \, \rho^{\,k}, \, k \in \n$.
How to identify all those ``fastest'' trajectories? 
\smallskip

\noindent \textbf{Problem 1}. {\em  How to 
characterise all switching laws~$A(\cdot)$ realising the maximal growth
of trajectories of the linear switching system? }

\smallskip

 We are going to show  that if the system has 
finitely many dominant products, then Problem~1 can be explicitly solved: 
 \smallskip 
 
{\em Suppose a switching system has 
a finite set of  dominant products and each of them has a unique and simple leading eigenvalue; then a switching law~$A(\cdot)$ generates trajectories of the maximal growth precisely when it is eventually periodic with the period equal to a dominant product. For all other laws,  
we have $\|A(k)\cdots A(1)\| \, = \, o(\rho^k)$ as $k \to \infty$}. 
\smallskip 

This condition  means that there exist numbers $n$ and $N$ 
and an infinite sequence of indices 
$s_1s_2\ldots $ such that $s_{k} = s_{k+n}$ for all $k > N$, 
$A(j) = A_{s_j}, \, j\in \n$,  
and the period $\Pi = A_{s_{N+n}}\cdots A_{s_{N+1}}$ is a dominant product.  

If we normalize the family so that  $\rho(\cA) = 1$, then Problem~1 becomes  to {\em characterise all switching laws that do not 
tend to zero} as $k\to \infty$. Since the normalized family has the same set of 
switching laws of the fastest growth, it suffices to consider the case~$\rho(\cA) = 1$. 

\begin{theorem}\label{th.120}
Let a system~$\cA$ be irreducible, normalized as~$\rho(\cA) = 1$, and  have  finitely many dominant products. Let also  
 the leading eigenvalues of all dominant products be unique and simple. 
Then 
all trajectories  of the system converge to zero apart from those corresponding to 
eventually periodic switching laws with a period 
equal to a dominant product.  
\end{theorem}

Thus, the switching laws of the maximal growth are precisely those eventually periodic ones 
with the period equal to a dominant product. All other switching laws tend to zero. 
\begin{ex}\label{ex.60}
{\em For the family~(\ref{eq.ex10}) from Example~\ref{ex.10},  
the dominant product is $A_1^3A_2$. Hence the trajectories 
of the fastest growth all have the form~$(A_1^3A_2)^j \Pi_0 \bx_0$, 
where $\Pi_0$ is an arbitrary product of the matrices~$A_1, A_2$. 
In particular, the trajectories 
$(A_1^3A_2)^jA_2\bx_0$ and $(A_1^3A_2)^jA_2A_1\bx_0, \, j \in \n$, are both of the fastest growth. The trajectories of the corresponding normalized family 
are shown in Fig.~\ref{fig60}. 
The left figure  presents  these two trajectories of the fastest growth:  brown and  green broken lines respectively. The points of all trajectories are shown in black. 
In Fig.~\ref{fig60} (right), one of remaining trajectories (of not the fastest growth) is shown.

\begin{figure}[ht!]
\begin{minipage}[h]{0.8\linewidth}
\center{\includegraphics[width=1\linewidth]{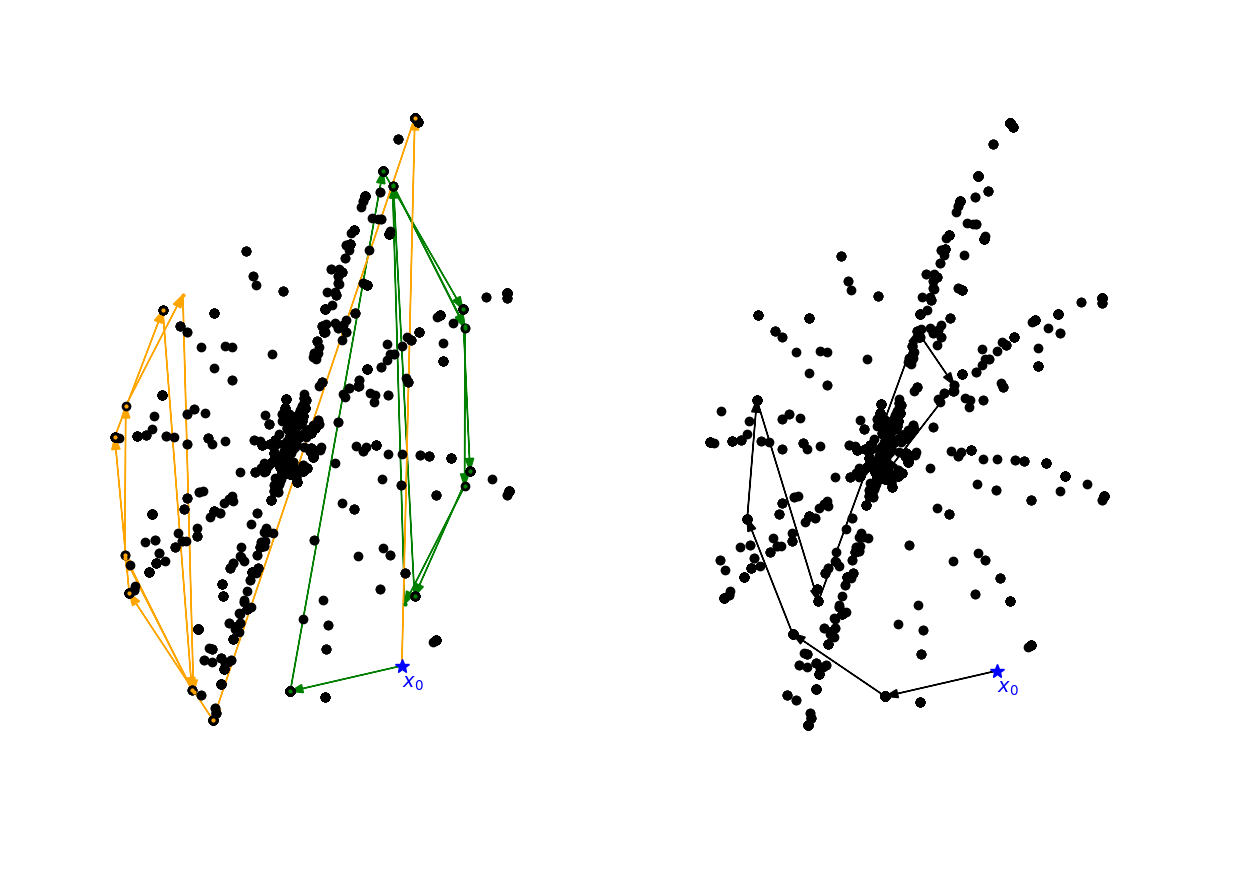}} 
 
\end{minipage}
\caption{{\footnotesize \textbf{Left:} Two 
fastest growth trajectories $(A_1^3A_2)^jA_1\bx_0\, $ and $\, (A_1^3A_2)^jA_0A_1\bx_0, \, j \in \n$, (brown and green respectively) for the family~(\ref{eq.ex10}); 
\textbf{Right:} some of other trajectories.}}
\label{fig60}
\end{figure}
}
\end{ex}

\begin{remark}\label{r.45}
{\em Theorem~\ref{th.120} characterises all switching laws 
for which $\|A(k)\cdots A(1)\| \ge C\, \rho^k$. Not all trajectories 
of those switching laws have the maximal growth. This property may depend on 
the initial point~$\bx_0$. Nevertheless, all trajectories of the 
maximal growth can be explicitly characterized. For the sake 
of simplicity, assume again that we have a normalized family, for which $\rho(\cA) = 1$. 
Since $\|\bx(k)\| \, \le \, \|A(k)\cdots A(1)\|\cdot \|\bx_0\|$, we see that 
if the switching law is not eventually periodic with a dominant period, then 
$\|\bx(k)\| \to 0$ as $k \to \infty$, hence this trajectory is not of the maximal growth. 
If it is eventually periodic with a dominant product~$\Pi$ of length~$|\Pi| = n$ as a period, then
the product $A(k)\cdots A(1)$  for $k = jn +N$ has the form 
$\Pi^j\Pi_0$, where $\Pi_0$ is a product of length~$N$. 
The leading eigenvalue of~$\Pi$ is equal to~$\rho(\cA) = 1$ and, by the assumption, this
eigenvalue is unique and simple. Denote by~$L$ the subspace of~$\re^b$ of dimension 
$d-1$ spanned by all vectors of the Jordan basis of~$\Pi$ except for the leading eigenvector. 
Actually,~$L$ is an orthogonal complement of the leading eigenvector  of the 
transpose matrix~$\Pi^T$. Then, if $\Pi_0\bx_0 \in L$, then 
$\Pi^j \Pi_0\bx_0 \to 0$ as $j\to \infty$, and hence 
$\bx(k) \to 0$ as $k \to \infty$. Otherwise, 
$\|\Pi^j \Pi_0\bx_0\| \ge  C$ for all $j$ and hence the norms $\|\bx(k)\|$
are bounded below by a positive constant for all~$k \in \n$. Thus, 
the trajectory~$\{\bx(k)\}_{k=0}^{\infty}$ has the maximal growth if and only if 
the switching law is periodic with a dominant period~$\Pi$ and 
$\Pi_0\bx_0 \notin L$. 
This gives the complete classification of all trajectories of the fastest growth. 
}
\end{remark}

Before giving a proof of Theorem~\ref{th.120} we need to introduce some more notation.  
\smallskip

{\em The cyclic tree of matrix products}. 
To a family of operators $\tilde \cA \, = \, \{\tilde A_1, \ldots
, \tilde A_m\}$ and to some product $\tilde \Pi \, = \, \tilde
A_{d_n}\cdots \tilde A_{d_1}$ with the spectral radius~$1$ we associate
the {\em cyclic tree} $\cT$ generated by the word $d_1\ldots d_n$
(or, which is the same, by the product~$\tilde \Pi$). 
It is defined as
follows. The root is formed by a cycle $\cR$ of $n$
nodes~$V_1, \ldots , V_n$. They are, by definition, the nodes of
zero level. For every $i\le n$ an edge (all edges are directed)
goes from $V_i$ to $V_{i+1}$, where we set $V_{n+1} = V_1$. At each
node of the root $m-1$ edges start to nodes of the first level.
So, there are $n(m-1)$ different nodes on the first level. The
sequel is by induction: there are $n(m-1)m^{k-1}$ nodes of the
$k$th level, $k\ge 1$, from each of them $m$ edges (``children'')
go to $m$ different nodes of the $(k+1)$st level.

Each index (letter) $q$ belongs to the alphabet $\{1,
\ldots , m\}$ and is associated to the matrix~$A_q \in \cA = \{A_1, \ldots , A_m\}$. Let us recall that we write products in the inverse order: from the right to the left. 
We assume that the root $\cR$ is primitive, i.e.,
is not a power of a shorter word. To every edge of the tree~$\cT$ we
associate a letter $q$ (or the corresponding matrix~$A_q$), as follows: the edge $V_iV_{i+1}$ of the root 
corresponds to $d_i\, , \, i=1, \ldots , n$; at each node~$V_i$ of the root $m-1$
edges start to the first level associated to all the~$m$ letters
except for~$i$. From each node $V$ of level~$k\ge 1$ exactly $m$
edges start associated to all the letters $1, \ldots , m$. 

We identify the words with the corresponding products of matrices
from~$\cA$. 
To a given point~$V_i \in \cR$ and to a given finite  word 
$q_1\ldots q_k$ we associate  the node $A_{q_k}\cdots A_{q_1}V_j$, which is the end of the
path from $V_i$ along the edges $q_1, q_2, \ldots , q_k$ respectively. 
The empty word corresponds to $V_i$. To an infinite word and to a node $V_i \in \cR$ we 
associate an infinite path~$V_i = V^{(0)}\to V^{(1)}\to V^{(2)}\to \cdots $ on the tree (all the paths are 
without backtracking) starting at~$V_i$. This path corresponds to the 
starting node~$V_i$ and to an infinite word $s_1 \ldots s_k \ldots $. 
A node $V^{(k)}$ on this path 
on $k$th level is $V^{(k)} \, = \, A_{s_k}\ldots A_{s_1}V_i$. 

The routine of Algorithm~1 can be described in terms of the tree~$\cT$.  
First, we have a root~$\cR = \{V_1, \ldots , V_n\}$. 
 At the first step
we take any node $V_i \in \cR$ and consider successively  its $(m-1)$
children from the first level. For each child  $\, V \, = \,
\tilde A V_i$, where $\tilde A \in \tilde \cA \setminus \{\tilde
A_{d_i}\}$ we  determine, whether or not 
$V$ belongs to the interior of~$G_1 = {\rm co}\, \cV_1$. If
it does, then $V$ is a {\em dead node} or {\em dead leaf} generating a {\em dead branch}:
we will never come back to $V$, nor to nodes of the branch
starting at $V$ (so, this branch is cut off). If it does not, then
$V$ is an {\em alive leaf}, and we add this element~$V$ to the set $\cV_1$
and to the set $\cH_1$. After the first iteration all alive
nodes of the first level form the set~$\cH_1$.  At the second
step we deal with the nodes from $\cH_1$  only and obtain the next set
of alive nodes of the second level~$\cH_2$, etc. Thus, after the
$k$th iteration we have a family $\cH_k$ of alive nodes from the $k$th
level, and a set $\cV_k \, = \, \cup_{j=0}^{k}\, \cH_j$. A node~$V$
belongs to the set~$\cV_k$ if and only if its level does not exceed~$k$ and it
belongs to an alive branch starting from the root. The convex body
$G_k$ is the convex hull ${\rm co} \, \cV_k$. The
convex body $G_{k-1}$ is invariant  if $\cH_{k}\, = \, \emptyset$,
i.e., the $k$th iteration produces no alive leafs (only dead ones).
This means that there are no alive paths of length~$k$ from
the root. Therefore $G_{k} = G_{k-1}$. Otherwise, if $\cH_{k}$ is nonempty, 
the algorithm makes  the next
iteration and goes to the $(k+1)$st level: we take children of each element
of $\cH_k$, determine whether they alive of dead, etc.
\smallskip 

Algorithm~2 works simultaneously with $r$ cyclic trees. Each 
cyclic tree~$\cT^{(i)}$ is generated by the $i$th candidate product~$\tilde \Pi^{(i)}, \, \ i = 1, \ldots, r$. On $k$th iteration we run over 
the set~$\cH_{k-1}$ that consists of nodes of the $(k-1)$st level of all the trees added in the $(k-1)$st iteration. The alive leafs of every node are added to~$\cH_k$, the dead leafs are omitted together with edges growing from them. 
When the whole set~$\cH_{k-1}$ is exhausted, we set $\cV_k = \cV_{k-1}\cup 
\cH_k$ and go to the next iteration.

\smallskip

{\tt Proof of Theorem~\ref{th.120}}. 
For the sake of simplicity, we assume that 
all the dominant products have real eigenvalues. Otherwise we replace 
the leading eigenvectors of some products by leading ellipses. 

By Theorem~\ref{th.25} (Section~6), for every family with finitely many dominant products, 
Algorithm~2 converges within finite time. Let it perform 
$k$ iterations. Denote by~$N$ the sum of the number $k$ and of maximal 
length of all dominant products of~$\cA$.  
 Take an arbitrary dominant product~$\Pi$ of length~$n$, denote the corresponding word by~$\pi$  and denote the leading eigenvector of~$\Pi$ by~$\bv$.  
Every infinite switching law which does not have period~$\pi$ has the 
form $\pi^{\ell}qs$, where $\ell \ge 0$ is an integer, $q$ is a word of length~$N$
whose prefix of length~$n$ is different from~$\pi$, and $s$ is an infinite word.
Let $Q \in \cA^N$ be the product corresponding to the word~$q$. 
We have $Q\Pi^{\ell}\bv \, = \, Q\bv$. By Algorithm~1, the trajectory 
starting at~$\bv$ with the switching law~$q$ has a point with G-norm 
(i.e., with the norm that has a unit ball~$G$) 
strictly less than one.  This point is a dead node on the path along the cyclic 
tree~$\cT$ generated by the product~$q$ and starting at~$\bv$. Denote by~$\mu$ the maximal 
G-norm of dead vertices of the trees generated by 
the dominant products. Since this set of vertices is finite, it follows that $\mu<1$. 
At every path starting from the root, a dead node has to appear 
by the~$k$th iteration, hence it corresponds to a product of length at 
most~$k + n \le N$. Since the  G-norm  is non-increasing along any trajectory, it follows that for every prefix~$q'$ of the word~$q$
of length at least~$N$, we have $\|Q'\Pi^{\ell}\bv\|_{G} \le \mu$. Thus, for every  switching law, unless it is eventually periodic with a dominant 
period,  there is a number $N_1$
such that for every its prefix~$s$ of lengths bigger than~$N_1$, we have 
 $\|S\bv\|_G \le \mu$. Choosing the maximum of those numbers~$N_1$ over all 
 vertices of~$G$ and taking into account that the G-norm of every linear operators  is achieved at one of the vertices of~$G$ we conclude that 
   for every switching law, unless it is eventually periodic with a dominant 
period,  every sufficiently long  its prefix~$s$ satisfies $\|S\bv\|_G \le \mu$. Therefore, every switching law, unless it is eventually periodic with a dominant 
period,   can be 
split into finite words such that the norms of the corresponding matrix products  are less than~$\mu$. This implies that 
the corresponding trajectory~$\{\bv_j\}_{j \in \n}$ contains a subsequence of points such that $\|\bv_{j_k}\|_{G}  \, \le \, \mu^k\|\bv\|$. 
Since the $G$-norm does not increase along any trajectory, it follows 
that $\|\bv_{j}\|_{G}  \, \le \, \mu^k\|\bv\|$ for all $j\ge j_k$. 
Hence,  $\bv_j$ tends to zero as $j\to \infty$.

{\hfill $\Box$}
\smallskip

\bigskip 

\begin{center}
\large{\textbf{8. The case of rational $\, {\rm mod}\, \pi\, $ argument}}
\end{center}
\bigskip

By the results of Sections 2 and 4, if a system has a dominant product, then 
its Barabanov norm is unique, provided the leading eigenvalue of the  dominant product is either real or complex with a rational $\, {\rm mod}\, \pi\, $ argument.
In the former case the unique norm is piecewise-linear, in the latter  case it is 
piecewise-quadratic. In both cases the norm has a simple structure and can be constructed with Algorithm~1. 
What can be said in the last case  when 
a dominant product has a non-real leading eigenvalue  with a rational $\, {\rm mod}\, \pi\, $ argument?  
We are going to see that in this case the set of Barabanov norms is much richer and 
more complicated. Nevertheless, we will classify all those norms (Theorem~\ref{th.40}).

First of all, the uniqueness may fail. For example, if  
$A_1$ is a rotation of the plane~$\re^2$  by $90^{\circ}$ and 
$A_2$ is an arbitrary operator with the (Euclidean) norm at most~$\frac12$, then the pair 
$\{A_1, A_2\}$ has infinitely many invariant bodies: every regular 
$4n$-gon is invariant, $n \in \n$. The operator~$A_1$ is dominant
and all other conditions of Theorem~\ref{th.10} (Section~2) are fulfilled, but the 
argument of the leading eigenvalue~$\lambda \, = \, e^{\frac{\pi i}{2}}$ is rational $\, {\rm mod}\, \pi$.  

In spite of the non-uniqueness, one may hope that all invariant bodies 
can still be constructed by~Algorithm~1 with a proper choice of the 
root sets~$\cR = \{V_j\}_{j=1}^n$. To introduce the idea we need one more notation  generalizing leading eigenvectors or leading ellipses. In this section we deal with 
solid ellipses $E = \Phi_{\bx, \by}(D)$, where $D\subset \re^2$ is a unit disc, 
and use for them the same notation as for the curves $E = \Phi_{\bx, \by}(\Gamma)$, where 
$\Gamma\subset \re^2$ is a unit circle. 

\begin{defi}\label{d.50}
For a given product $\Pi \in \cA^{\n}$ with a simple complex leading eigenvalue~$\lambda$, 
a convex compact subset~$V$ of the leading eigenspace of~$\Pi$ is called 
an admissible set if $V \ne \{0\}, \, V = -V$,  and $\Pi V \, = \, 
|\lambda|\, V$. 
\end{defi}
If $\lambda \in \re$, 
then an admissible set is a segment parallel to the leading eigenvector. 
If $\lambda$ is complex with an irrational $\, {\rm mod}\, \pi\, $ argument, then 
an admissible set is a leading ellipse $E = \Phi_{\bx, \by}(D)$, 
where $\bx + i\by$ is the leading eigenvector and $D$ is a unit disc on the plane.
If the argument of~$\lambda$ is rational $\, {\rm mod}\, \pi$, then there are 
infinitely many, up to homothety, admissible sets. The proof of the following lemma is omitted since it is simple.  
\begin{lemma}\label{l.20}
Let a product $\Pi \in \cA^{\n}$ has a complex leading  eigenvalue
with an incommensurable with $\pi$  argument~$\varphi$; then 
a  set~$V$ is  admissible if and only if 
$V = \Phi_{\bx, \by}(M)$, where $M \subset \re^2$ is a 
convex body symmetric about the origin and mapped to itself 
by the rotation by the angle~$\varphi$.   
\end{lemma}
The generalization of Algorithm~1 to an arbitrary 
admissible starting set~$V_1 = \Phi_{\bx, \by}(M)$ is  the following.  We 
take an admissible set~$V_1$ and define as usual the sets 
$V_j = \tilde A_{d_{j-1}} \cdots \tilde A_{d_1}V_1, \, j = 2, \ldots , n$. 
Since the restriction of the operator~$\Pi = A_{d_n}\cdots A_{d_1}$ is a composition of 
 rotation by the angle~$\varphi$ and of multiplication 
 by~$|\lambda| = \nu(\Pi)$, we see that the normalized operator $\tilde \Pi$ is a rotation by the angle~$\varphi$. Hence $\tilde \Pi V_1\, = \, V_1$, so the 
 sets~$\{V_j\}_{j=1}^n$ 
 indeed form a cycle~$V_1\to \cdots \to V_n \to V_1$ with the edges 
  (operators)~$\tilde A_{d_1}, \ldots , \tilde A_{d_n}$. Thus, Algorithm~1 
  with the candidate product~$\Pi$ and with the root~$\cR = \{V_j\}_{j=1}^n$ 
  produces an invariant 
  body~$G_k$ whenever it halts after $k$th iteration. 
  
  Thus, taking an arbitrary admissible set~$V_1$ we define the root $\cR$ and 
  start Algorithm~1. If $M=D$ is  a disc,  then $V_1$ is the leading ellipsoid. 
  In this case, as it follows from    Theorem~A,   Algorithm~1 halts within finite time provided the product~$\Pi$ is dominant.   
  However, this is not true for some other admissible sets~$V_1$ as  Example~\ref{ex.100} below demonstrates. Therefore, this direct generalization of Algorithm~1 may not be applicable
  for admissible sets other than ellipses.

\begin{ex}\label{ex.100}{\em 
We are going to construct  a pair of $4\times 4$ matrices with one dominant 
product and a complex leading eigenvalue with a rational $\, {\rm mod}\, \pi\, $ argument
and an admissible  set~$V_1$ 
for which Algorithm~1 does not terminate within finite time. 

\smallskip 

We consider the space~$\re^4$ and its 
two-dimensional subspace~$L\, = \, \{(x_1, x_2, 0, 0)^T \in \re^4\}$. Sometimes we denote a point from~$L$
as $\bx = (x_1, x_2)^T$.  An orthogonal projection of a point~$\by$ to~$L$
is denoted by~$\tilde \by$. 
Consider a regular hexagon~$H$ on~$L$ 
centered at the origin with one vertex at the point~$\bv = (0,1,0,0)^T$. Its side is 
equal to one. 
 Take small~$\tau > 0$ and a vector $\bb \, = \, \bigl(\tau, 1-4\tau^2, \tau, \tau \bigr)^T \in \re^4$. For all sufficiently small~$\tau$, we have 
 $\|\bb\|_2 < 1$. 
 Consider a pair of matrices $\cA = \{A_1, A_2\}$ with 
 \begin{equation}\label{eq.hex}
 A_1 \ = \ 
 \left( 
 \begin{array}{cccc}
 -\frac{1}{2}& -\frac{\sqrt{3}}{2}& 0 & 0\\
  \frac{\sqrt{3}}{2}&-\frac{1}{2}&  0 & 0\\
  0& 0& \frac12 & 0\\
  0 & 0& 0 & \frac14 
 \end{array}
 \right)\ ; 
 \qquad 
 A_2 \ = \ \bb \, \bb^T\, . 
 \end{equation}
 Thus, the matrix $A_1$ consists of three diagonal blocks. The 
 first $2\times 2$ block (we call it~$B$) is the rotation of the plane~$L$ by $120^{\circ}$, the other two one-dimensional blocks are $\frac12$ and $\frac14$. The rank-one matrix $A_2$ defines the operator $A_2\bx  = (\bx , \bb) \, \bb$. This is an orthogonal projection to the direction of vector~$\bb$ multiplied by 
 $\|\bb\|^2$. Clearly, $\|A_2\| = \|\bb\|^2 < 1$. On the other hand, 
 $\|A_1\| = \rho(A_1)= 1$. Hence, for every product~$S$ of matrices $A_1, A_2$, 
 we have $\|S\| \le \|\bb\|^2$, unless $S$ is a power of~$A_1$. Therefore, $A_1$
 is a (unique!) dominant product of the family $\cA = \{A_1, A_2\}$. 
 Consequently, $\cA$ has a spectral gap~$\bigl( \|\bb\|^2 , 1\bigr)$ 
 and $\rho(\cA) = \rho(A_1) = 1$. }
 \end{ex}
\begin{prop}\label{p.20}
For every small~$\tau \, > \, 0$, Algorithm~1 applied to the
 pair~$\{A_1, A_2\}$ with the initial admissible  set~$V_1 = H$, makes infinitely many iterations and produces 
an invariant body with infinite discrete set of extreme points. 
\end{prop}
\begin{remark}\label{r.20}
{\em Algorithm~1 applied to the pair~$\{\bv, - \bv\}$ gives the same result as 
being applied to the hexagon~$H$. 
}
\end{remark}
{\tt Proof of Proposition~\ref{p.20}}.  For every $k$, the points  $\pm \, A_1^k\bv$
are vertices of~$H$. Denote $\bc = A_2\bv\, = \, (1-4\tau^2)\bb$. 
For small $\tau$, the point $\tilde \bc \, = \, 
\Bigl((1-4\tau^2)\tau, (1-4\tau^2)^2\Bigr)^T$ (the projection of~$\bc$ to~$L$)
is out of~$H$. 
Denote by $H_{\tau}$ the regular hexagon in~$L$ centered at the origin 
and having  one of the  vertices at the point~$\tilde \bc$. All vertices of $H_{\tau}$ 
are out of~$H$. 
Clearly, all points $\, \pm \, B^k \tilde \bc, \, k \in \n$, are also 
vertices of~$H_{\tau}$. 
Finally, if $\tau$ is small enough, then the point~$\bv$
has the biggest in modulus scalar product with the vector~$\bb$ among all vertices of 
the hexagons~$H$ and $H_{\tau}$ and the vector~$\bb$.  
So, the maximum of the functional $F(\bx) = (\bb, \bx)$ on the set 
${\rm co}\, \{H, H_{\tau}, \bb\}$ is attained at a unique point~$\bv$. 
Therefore, projections of all points generated by~Algorithm~1 to the 
plane~$L$
are in the set~${\rm co}\{H, H_{\tau}\}$.  
Consider the sequence
$$
A_1^{3k}A_2\bv \ = \ \bigl( B^{3k}\tilde \bc \, , \, 2^{-k}\tau \, , \, 4^{-k}\tau  \bigr)^T\ =  \ \bigl(\tilde \bc \, , \, 2^{-k}\tau \, , \, 4^{-k}\tau  \bigr)^T\, , 
\quad k \ge 0 \, . 
$$
All these points are convex independent (none of them is in the convex hull of others) since so are the points $\bigl( 2^{-k}\tau \, , \, 4^{-k}\tau  \bigr)^T, \ 
k \in \n$, 
because they all lie on the positive part of the parabola 
$y \, = \, \frac{1}{\tau}\, x^2,  \, x \ge 0$. If some point $A_1^{3j}A_2\bv$ is not an extreme point of 
the body~$G$ generated by the algorithm, then by the Minkowski theorem it must be a convex combination of other extreme points. However, the projection of~$A_1^{3j}A_2\bv$ 
to $L$, which is the point~$\tilde \bc$, is extreme for the projection 
of the set~$G$ to~$L$, which is the set 
${\rm co}\, \{H, H_{\tau}\}$. Hence,~$\tilde \bc$ must be a convex combination of 
points generated by Algorithm~1 whose projection to $L$ coincides with~$\tilde \bc$, i.e., points from the sequence  $\{A_1^{3k}A_2\bv\}_{k \ge 0}$. 
This is impossible due to convex independence of this sequence. 

Thus, Algorithm~1 starting with the set~$V_1 = \cH$ produces a sequence 
of extreme points $\{A_1^{3k}A_2\bv\}_{k \in \n}$. Other sequences are 
$\{A_1^{3k+1}A_2\bv\}_{k \ge 0}, \, \{A_1^{3k+2}A_2\bv\}_{k \ge 0}$
and the sequences symmetric to them about the origin. 
Those six sequences converge to vertices of the hexagon~$H_{\tau}$. 
The convex hull of these six sequences and of vertices of $H$ and of $H_{\tau}$
is the invariant body~$G$ produced by Algorithm~1. This body has an infinite discrete set of extreme points.

{\hfill $\Box$}
\smallskip

Thus, a direct application of  Algorithm~1 to an arbitrary admissible 
set~$V_1$ may lead to divergence. Therefore, a classification of invariant sets 
in the case of complex leading eigenvalue with a rational $\, {\rm mod}\, \pi\, $ argument requires a different
procedure. This can be done by modifying Algorithm~1 as stated below. 
We describe the modified algorithm as Algorithm~3. Each iteration of the new algorithm deals with infinite sets of points, therefore it cannot be considered as 
a finite procedure and its significance is rather theoretical. Nevertheless, it establishes a complete classification of 
invariant sets and of Barabanov's norms in the case of rational $\, {\rm mod}\, \pi\, $ argument.  
It is realised in the same way as Algorithms~1 but with two 
differences:

\smallskip 

1) The starting set~$V_1$ is an arbitrary admissible set for the 
candidate product~$\Pi$. In particular, for $V_1 = E_1$, 
we obtain 
Algorithm~1 in case of complex leading eigenvalue.  
\smallskip 

2) Every node~$V$ of the cyclic tree~$\cT$ is either 
an element of the root~$V_j \in \cR\, = \, \{V_1, \ldots , V_n\}$ or 
the end of a finite path starting at some node~$V_j \in \cR$. 
Denote by $\tilde \Pi_j$  the $j$th cyclic permutation of $\tilde \Pi$, 
which sends $V_j$ to itself. Let $\pi_j$ be the word corresponding to the product~$\tilde \Pi_j$ and $\pi_j^{\infty} = \pi_j \pi_j \ldots $ be the corresponding infinite word. 
For a node~$V \in \cT$, we denote by $\pi_j^{\infty}(V)$ the corresponding infinite path 
$\pi_j^{\infty}$ along~$\cT$ starting at the node~$V$. 
\smallskip 

In one step we add the following sets to~$\cH_{k+1}$: 
\smallskip 

a) all nodes of the infinite path $\pi_j^{\infty}(V)$; 
\smallskip 

2)  the $m-1$ children $\tilde A_kV, \, k \in \{1, \ldots , m\}, \, \, A_k$
is not the first matrix in the product~$\Pi_j$ (i.e., 
the child does not belong to this path). 
\smallskip 

So, each step adds infinitely many nodes. 
We check all those new nodes. 
A node~$V'$ is dead 
if and only if it is (all its points) is in the interior of the  current set~$G_k$. 
In contrast to Algorithm~1, here the set of new vertices $\cH_k$ may contain an infinite set of nodes and is not necessarily located in one level. 
\smallskip

Now write the formal routine.  
\bigskip

\textbf{Algorithm 3}. 
\smallskip 

\textbf{I. Choosing the candidate product}. The same as in~Algorithm~1.  
 \smallskip 
 
 \textbf{II. The routine.} 
  \smallskip 
  
Choose $V_1$ is an arbitrary subset of leading eigenspace~$L_1$ of $\tilde \Pi$ such that $\tilde \Pi V_1 = V_1$ and $V_1$ is symmetric about the origin. 
Then define the root 
 $\cR = \{V_1, \ldots , V_n\}$ from the set~$V_1$ as in Algorithm~1. 
Then we construct a sequence of  sets $\cV_i$ of nodes and their
subsets~$\cH_i \subset \cV_i$ (may be infinite) as follows:
 \smallskip 
 
{\tt Zero iteration}. We set $\cV_0 = \cH_0 = \cR$.
 \smallskip 
 
{\tt $k$th iteration, $k \ge 1$}. We have a  set of nodes $\cV_{k-1}$ and its subset
$\cH_{k-1}$.
We set $\cV_{k} = \cV_{k-1},  \, \cH_k = \emptyset$. 
Take an arbitrary node~$V \in \cH_{k-1}$. It is the end of a finite path 
starting at a node~$V_j$ of the root. 
Denote by $\tilde \Pi_j$  the $j$th cyclic permutation of $\tilde \Pi$, 
which sends $V_j$ to itself. For every $\tilde A \in \tilde \cA $, 
which is different from the first matrix of the product~$\tilde \Pi$, 
check whether $ \tilde A\, V$ is in the interior of
$G_{k-1} \, = \, {\rm co} \{V \ | \ V \in \cV_{k-1}\}$. 
If it is, then we omit the set $ \tilde A\, V$ and take the next pair $(V , \tilde \cA) \in \cH_{k-1}\times  \tilde \cA$, otherwise we add $ \tilde A \, V$ to $\cV_{k}$ and to $\cH_k$. If $k\ge 2$, we do this for all pairs from~$\cH_{k-1}\times  \tilde \cA$, except for those where $A$ is the first element of~$\tilde \Pi_j$. 
 If $k = 1$ and hence $\cH_{k-1} = \cH_0 = \cR$ and $V=V_j \in \cR$, then we exclude $n$ pairs~$(V_j, \tilde A_j), \, j = 1, \ldots , n$. Finally, 
 if $A$ is the first matrix of the product~$\tilde P_j$, then 
 we consider the infinite path~$\pi_j^{\infty}(V)$. 
 Take the highest (i.e., on the maximal level) node~$V'$ of this path which belongs to ${\rm int}\, G_k$. 
 We remove this node and the corresponding branch of the tree growing from it, 
 including the remainder of this path. All the nodes of this path higher than~$V'$
 are added to both $\cV_k$ and to $\cH_k$. If such a node~$V'$ does not exist, then 
 all nodes of~$\pi_j^{\infty}(V)$ are added to $\cV_k$ and to $\cH_k$.

When all proper pairs $(V,  \tilde A)$ are exhausted, both $\cV_{k}$ and $\cH_k$ are constructed.
We define  $G_k = {\rm co} \{V \ | \ V \in \cV_{k}\}$ and  have

{\tt Termination}. The algorithm halts when $\cV_{k} = \cV_{k-1}$, i.e., $\cH_k = \emptyset$. In this case $G_k$  is an invariant convex body
for~$\cA$.
\smallskip 

\hfill \textbf{End of the algorithm}.

\begin{remark}\label{r.50}
{\em Algorithm~3 is rather theoretical because each iteration 
 assumes infinite number of steps: verifying the assertion~$\tilde A V \in {\rm int}\, G_{k-1}$ for infinitely many nodes~$V \in \cH_{k-1}$. Nevertheless, it shows the 
 theoretical way to find the invariant convex body generated by an arbitrary admissible
 set~$V_1$. On the other hand, if $V_1 = E_1$ (the leading ellipsoid), then Algorithm~1 converges within finite time provided~$\Pi$ is dominant. Therefore, in this case there is no need to 
 apply Algorithm~3. Moreover, it is not reasonable to apply Algorithm~3 for 
 computing the joint spectral radius either, because the JSR can always be computed with Algorithm~1 for~$V_1 = E_1$.

}
\end{remark}
\smallskip

\begin{theorem}\label{th.30}
Let a family $\cA$ possess a unique dominant product~$\Pi$ and let 
$\Pi$ have a unique and simple complex leading eigenvalue 
with a rational $\, {\rm mod}\, \pi\, $ argument. Then for every admissible set~$V_1$, 
Algorithm~3 terminates  within finite number of iterations and produces an invariant body. 
\end{theorem}
{\tt Proof}. If the algorithm does not converge within finitely many iterations, 
then there is an infinite path on the cyclic tree~$\cT$
starting at the root that consists of alive nodes and is constructed by infinitely 
many iterations. 
Denote the  node of this path on $k$th level by~$V^{(k)}$. 
Since $V^{(0)}$ belongs to the root,  it can be assumed  that $V^{(0)} = V_1$. 
The set~$V^{(k)}$, which is the image of~$V_1$ by the corresponding matrix product of length~$k$, is alive if it is not in the interior of the body~$G_k$ 
constructed in the $k$th iteration. Therefore, the 
diameter of the set~$V^{(k)}$ cannot converge to zero as~$k\to \infty$.  In view of Theorem~\ref{th.120} from Section~7, this means that the path $V^{(0)} \to V^{(1)} \to \cdots $
corresponds to an eventually periodic switching law with the period~$\pi$ (the word associated 
to the product~$\Pi$). 
Suppose the periodic part starts after $j$th iteration, at the node~$V^{(j)}$. Then 
this is the infinite path~$\pi^{\infty}(V^{(j)})$. However, all nodes of this path are added at once 
 in the $j$th iterations. Thus, the whole path $V^{(0)} \to V^{(1)} \to \cdots $
 is constructed in the first $j$ iterations, which contradicts to the assumption.

{\hfill $\Box$}
\smallskip

Thus, for every admissible set~$V_1$,  Algorithm~3 converges and produces an invariant body. 
The next result shows that every invariant body is obtained this way. 

\begin{theorem}\label{th.40}
Let a family $\cA$ possess a unique dominant product~$\Pi$ whose leading 
 eigenvalue has a rational $\, {\rm mod}\, \pi\, $ argument. 
Then every its invariant body is constructed by Algorithm~3 with some 
 admissible subset~$V_1$ of the leading plane. 
\end{theorem}
{\tt Proof}. 
Assume after possible normalization that $\rho(\Pi) = 1$. 
Let~$G$ be an invariant body and $V$ be its intersection  with 
the leading eigenspace~$L$ of~$\Pi$. Then $\Pi V \subset V$. 
On the other hand, since the two-dimensional restriction of the operator~$\Pi$ 
to~$L$ has both its eigenvalues equal to one in modulus, 
it preserves the two-dimensional volume.  Hence~$\Pi V = V$ and so $V$ is admissible. 
Since both eigenvalues of~$\Pi|_{V}$ have rational $\, {\rm mod}\, \pi\, $l arguments, it follows that 
each point~$\bx \in \partial V$ is recurrent for the family~$\cA$. On the other hand, by  Proposition~\ref{p.10} from~Section~5, all recurrent points are on~$\partial V$. Therefore,~$\partial V$ 
is the locus of recurrent points. Invoking now Theorem~B (Section~5) we conclude that 
$G$ is the closure of convex hulls of all trajectories starting from~$\partial V$. 
Hence~$G$ it is obtained  by Algorithm~3
from the candidate  product~$\Pi$ and the admissible set~$V_1 = V$.

{\hfill $\Box$}
\smallskip

\begin{remark}\label{r.60}
{\em Note that if we need one invariant set/Barabanov norm,  Algorithm~3 is not necessary. 
This can be done by Algorithm~1 with the admissible set~$V_1$
being the leading ellipse~$E_1$. Algorithm~3 is needed only to obtain all invariant sets/Barabanov norms.  

Theorem~\ref{th.40} classifies all invariant sets of the family~$\cA$ and explains why it may not be unique: every admissible set~$V_1$ generates an invariant body. 
If the leading eigenvalue  of~$\Pi$ has an irrational $\, {\rm mod}\, \pi\, $ argument, then there is a unique 
(up to multiplication by a constant) admissible set, which is the leading ellipse~$E_1$. 
In the case of rational $\, {\rm mod}\, \pi\, $ arguments there are many admissible sets. 

The transfer to Barabanov's norm is realised in the standard way: we take 
an arbitrary admissible set~$V_1^*$ for the dual family~$\cA^*$ and generate 
an invariant body~$G^*$ applying Algorithm~3. Then the Barabanov norm is~$f(\bx) \, = \, 
\max_{\by^* \in G^*}(\bx, \by^*)$.   
}
\end{remark}

\begin{remark}\label{r.70}
{\em If a family~$\cA$ has several dominant products, then 
if all of them have leading eigenvalues which are either 
real or complex with irrational $\, {\rm mod}\, \pi\, $ arguments, then 
all invariant bodies of~$\cA$ are convex hulls of finitely many 
points and ellipses (Corollary~\ref{c.10}). If 
at least one of dominant products, say, $\Pi^{(j)}$ has a non-real
 leading eigenvalue with a rational $\, {\rm mod}\, \pi\, $ argument, then there are more complicated 
invariant bodies. Namely,  the corresponding root $\cR^{(j)} = \{V_1^{(j)}, \ldots , V_{n_j}^{(j)}\}$ can be  generated by an arbitrary admissible set~$V_1^{(j)}$
of the product~$\Pi^{(j)}$.    
}
\end{remark}

\bigskip

\begin{center}
\large{\textbf{9. Barabanov norms for positive systems}}
\end{center}
\bigskip

A linear switching  system is called {\em positive} if 
all matrices of the family $\cA$ are (entrywise) non-negative.
If a positive system starts at a non-negative point~$\bx_0 \in \re^d_+$, 
then the whole trajectory is  in~$\re^d_+$.    
For positive  systems, the invariant polytope algorithm is extremely efficient even in dimension of $5000$ and higher~\cite{GP2, M}. However, to reach this efficiency we need 
to modify the concepts of invariant body and of Barabanov norm. We recall that inequalities $\bx \ge \nul, \bx \ge \by, 
A \ge \nul, A\ge B$
are understood entrywise. The positive orthant is $\re^d_+ \, = \, \{\bx \in \re^d \ |
\ \bx \ge \nul\}$. For positive systems, we usually work only with norms  defined in~$\re^d_+$. Moreover, it suffices to consider only {\em monotone} norms~$f$ for which 
$f(\bx) \ge f(\by)$ whenever $\bx \ge \by \ge \nul$. 
Respectively, we can consider {\em monotone convex bodies}~$G$
which lie in~$\re^d_+$ and possess the following property: 
if $\bx \in G$, then $\by \in G$ whenever~$\by \le \bx$. 
Similarly one defines the monotone convex hull of a set $K \subset \re^d_+$: 
$$
{\rm co}_{-}\,K\quad = \quad \Bigl\{\, \by \in \re^d_+ \ \Bigl|  \ \exists 
\ \bx \in {\rm co}\, K\, , \ \bx \ge \by  \, \Bigl\}. 
$$
Thus, the monotone convex hull contains  the usual convex hull plus
all points majorated by it.  A monotone convex hull of a 
finite set is a {\em monotone polytope}. 
In contrast to the usual polytope, a monotone polytope can have less than~$d$
vertices. For example, it  can have only one vertex~$\ba$, in which case it is a parallelepiped~$\{\bx \in \re^d \ | \ \nul \le \bx \le \ba \}$. 

A monotone norm is Barabanov if $\lambda f(\bx)\, = \, \max_{A_j \in \cA}f(A_j\bx)$
for all~$\bx \in \re^d_+$.  
A monotone convex body~$G$ is invariant for~$\cA$ if $\lambda G \, = \, 
{\rm co}_{-}\, \Bigl\{\, \cup_{A_j \in \cA} A_jG\, \Bigr\}$. 
The monotone invariant body and the monotone 
invariant norm are related by the {\em monotone polar transform}. 
The monotone polar to a set~$G \subset \re^d_+$ is 
$$
G^*_{-} \ = \ \bigl\{\, \bx \in \re^{d}_+ \ \bigr| \ 
 \sup_{\by \in G}(\bx , \by) \, \le \, 1  \, \bigr\} \, . 
$$
Note that for $\by \ge 0$ the relation $\bx_1 \le \bx_2$ implies that 
$(\bx_1 , \by) \le (\bx_2 , \by)$. Therefore, 
the sets~$G$ and ${\rm co}_{+} \, G $ have the same monotone polar. 
If $f$ is a monotone invarinat norm for~$\cA$, then the monotone polar to its unit ball is a monotone  invariant body
for~$\cA^*$~\cite{GP1}.

Finally, the irreducibility assumption for positive systems is weakened to {\em positive irreducibility}: the matrices from~$\cA$ do not share an invariant coordinate subspace i.e., subspace of the 
form $L_{S} \, = \, \{\bx \in \re^d \ | \ x_i = 0,  i \notin S\}$, where 
$S \subsetneq \{1, \ldots d\}$.
\smallskip 

\noindent \textbf{Theorem~D}~\cite{GP1}. {\em A positively irreducible system~$\cA$ possesses a monotone Barabanov norm and a monotone  invariant body. The unit ball of the 
monotone Barabanov norm is a polar to the invariant body of the transpose system~$\cA^*$.}
\smallskip 

 As for the structure of invariant bodies, the Perron-Frobenius theorem
 reduces the three possible cases of leading eigenvalues (real, complex with an irrational $\, {\rm mod}\, \pi\, $
 argument,  and non-real with a rational $\, {\rm mod}\, \pi\, $
 argument)  to one case. 
 Indeed, since a non-negative matrix always has a non-negative leading eigenvalue, 
 the cases of complex leading eigenvalues become impossible. Hence, Theorems~\ref{th.10} and~\ref{th.12} from Section~2 get the following simple form:
\begin{cor}\label{c.50}
If  a family of non-negative matrices~$\cA$ has a unique 
 dominant product with a unique and simple leading eigenvalue, then  it possesses a unique invariant body  and a unique Barabanov norm. 
  \end{cor}
Note that this unique invariant body 
may not be monotone. For the corresponding example, see, for instance~\cite[Figure~4]{GZ08}. 
However, a monotone invariant body does exist.    
  \begin{theorem}\label{th.50}
If  a family of non-negative matrices~$\cA$ has a unique 
   dominant product with a unique and simple leading eigenvalue, then  it possesses a unique monotone invariant body
   and a unique monotone Barabanov norm. The invariant body is a monotone polytope. 
The monotone Barabanov norm 
is given by the formula
\begin{equation}\label{eq.lin-p}
 f(\bx) \quad = \quad \max_{\bv^*} \, (\bv^*\, , \, \bx),\qquad \bx \in \re^d_+, 
 \end{equation}
 where the maximum 
 is taken over all vertices~$\bv^*$ of the monotone invariant polytope~$G^*$ of the dual family~$\cA^*$.  
\end{theorem}
The algorithm for construction of the monotone invariant polytope works in the same way as Algorithm~1 with the only difference: for each $k$, the polytope $G_k$ is a monotone convex 
hull of~$\cV_k$ (not just a convex hull as in~Algorithm~1). The proof of 
Theorem~\ref{th.50} is realized in the same way as for Theorem~\ref{th.10}. 
We only remark that if $G$ is the (usual) invariant body, then the 
monotone invariant body is the monotone convex hull of the set $G\cap \re^d_+$. 
Since a set and its monotone convex hull have the same monotone polar, it follows that 
the Barabanov norm restricted to~$\re^d_+$ is monotone. In particular, the 
Barabanov norm on~$\re^d_+$ coincides with the monotone Barabanov norm. 

As a rule, a monotone invariant polytope has much less vertices. 
In practice, even in very high dimensions, the number of vertices of  an invariant monotone polytope do not exceed several dozens. 
That is why in dimensions of several thousands the algorithm constructs 
the invariant body and the Barabanov norm within a few iterations. We report the numerical results  in the next section, Table~\ref{tab.20}.

{\hfill $\Box$}
\smallskip

\bigskip

\begin{center}
\large{\textbf{10. Numerical results}}
\end{center}
\bigskip

We report the results of performing Algorithms 1 and 2 for randomly generated matrices. 
Many results for matrices taken from practical applications can be found in~\cite{GP1, GP2, M}
and they are either similar or better than those for random matrices.   The numerical results presented  here are done 
by the most recent version of Algorithms 1 and 2 from~\cite{M}. 
Table~\ref{tab.10} shows the results of Algorithm 1 for arbitrary matrices with the case $\lambda \in \re$ (real leading eigenvalue of the dominant product). 
For even dimensions $d$ from $2$ to $20$, we took pairs of random $d\times d$ matrices
$\cA = \{A_1, A_2\}$ 
and normalise them either as  $\|A_1\| = \|A_2\|$ (the first column) or as $\rho(A_1) = \rho(A_2)$ (the second column). This normalization makes the problem more complicated, 
otherwise in most cases  the dominant product has length~$1$, i.e., either $A_1$ dominates $A_2$ or vice versa. 
For every dimension $d$ in each case $20$ experiments have been made in a standard laptop
and the median value of the computer time and of the number of vertices 
of the invariant polytope~$G$ is put in the table. The symbol $\hphantom{m}\#V$ denotes   
the number of pairs of vertices, so the invariant polytope has twice as many vertices. 
 Note that we did not remove the redundand verices, so the real number of vertices is 
usually much smaller. In the case  $\lambda \notin \re$ (complex leading eigenvalue of the dominant product), Algorithm~1 works slower. In our experiments it is mostly applicable 
for dimensions $\le 13$, for higher dimensions, the running time  often
exceeds reasonable limits. This can be explained by the fact that 
the conic programming problem~(\ref{eq.cp}) takes more time than the 
linear programming problem~(\ref{eq.cp}) in the real case. The total number or vertices (ellipses in this case) in the invariant body slightly exceeds the number if vertices
in  Table~\ref{tab.10} for the real case. 

\begin{table}[ht]
	\centering
	\caption{Computation of the Barabanov norm, arbitrary matrices}    
	\begin{tabular}{|c|rr|rr|}
		\hline
		&\multicolumn{2}{|c|}{{\footnotesize $\|A_1\| = \|A_2\|$}}&\multicolumn{2}{|c|}{{\footnotesize $\rho(A_1) = \rho(A_2)$}}\\
		{\footnotesize dim}&{ \footnotesize \hphantom{m}time}&{\footnotesize \hphantom{m}\#V}&{\footnotesize \hphantom{m}time}&{\footnotesize \hphantom{m}\#V}\\ \hline  
		2    &      1.1$\,s$ &     $5\cdot 2$        &      1.2$\,s$ &        $6\cdot 2$ \\   
		     
		4    &      1.4$\,s$ &     $17\cdot 2$      &      1.8$\,s$ &        $77\cdot 2$          \\  
		      
		6   &      2.0$\,s$ &     $47\cdot 2$      &      2.5$\,s$ &        $130\cdot 2$    \\ 
		       
		8   &      2.5$\,s$ &     $100 \cdot 2$       &      3.9$\,s$ &        $220\cdot 2$   \\  
		
		10   &      4.9$\,s$ &     $270\cdot 2$       &      5.1$\,s$ &        $320\cdot 2$  \\  
		      
		12   &      4.7$\,s$ &     $280 \cdot 2$       &      11$\,s$ &       $770 \cdot 2$    \\  
		      
		14  &      8.4$\,s$ &     $510 \cdot 2$     &       21$\,s$ &      $1100\cdot 2$  \\  
		      
		16  &       25$\,s$ &     $1100\cdot 2$      &       33$\,s$ &      $1400 \cdot 2$   \\ 
		
		18  &       90$\,s$ &     $2100 \cdot 2$       &       200$\,s$ &    $2500\cdot 2$   \\

		20  &       295$\,s$ &     $3100 \cdot 2$       &       5000$\,s$ &   $6200  \cdot 2$ \\ 
		
		       \hline
	\end{tabular}
	\label{tab.10}  
\end{table}

Table~\ref{tab.20} shows the results for non-negative matrices. In the first column the matrices are positive and in the second they are sparse with $90\%$ zero entries. We see that 
in the non-negative  case the algorithm is extremely efficient. Usually it constructs the Barabanov norm within $3-4$ iterations and this seems not to depend on the dimension.  
The number of vertices is usually around $8$ since we did not remove redundant vertices. 
For every  dimension $d$, in each case $20$ experiments have been made 
and the median values are reported. The algorithm always halted within finite time.

\begin{table}[ht]
	\centering
	\caption{Computation of the monotone Barabanov norm, non-negative matrices}    
	\begin{tabular}{|c|rr|rr|}
		\hline
		&\multicolumn{2}{|c|}{{\footnotesize 0\%~sparsity}}&\multicolumn{2}{|c|}{{\footnotesize 90\%~sparsity}}\\
		{\footnotesize dim}&{ \footnotesize \hphantom{m}time}&{\footnotesize \hphantom{m}\#V}&{\footnotesize \hphantom{m}time}&{\footnotesize \hphantom{m}\#V}\\ \hline  
		20    &      0.3$\,s$ &     7        &      1.7$\,s$ &        42           \\   
		     
		50    &      0.3$\,s$ &     8        &      1.6$\,s$ &        50          \\  
		      
		100   &      0.4$\,s$ &     8        &      0.8$\,s$ &        25     \\ 
		       
		200   &      0.5$\,s$ &     8        &      1.0$\,s$ &        23   \\  
		      
		500   &      1.2$\,s$ &     8        &      1.8$\,s$ &        16    \\  
		      
		1000  &      6.3$\,s$ &     8        &       11$\,s$ &        16    \\  
		      
		2000  &       35$\,s$ &     8        &       72$\,s$ &        16    \\        \hline
	\end{tabular}
	\label{tab.20}  
\end{table}

We see that for arbitrary matrices, the construction of Barabanov's norm in dimensions less than $15$ takes 
more or less the same time as for constructing other Lyapunov functions by known methods, which give only approximate values of JSR. For positive systems, Barabanov's norm is constructed much faster  even for very large dimensions. 
\bigskip 

 \textbf{Acknowledgements}. The author is grateful to T.Zaitseva to T.Mejstrik for their  help in making pictures and  for useful discussions of the computational issue.

\end{document}